\documentclass[11pt,reqno]{amsart}
\usepackage{amscd,amssymb,amsmath,amsthm}
\usepackage{mathtools}
\usepackage[english]{babel}
\usepackage[T1]{fontenc}
\usepackage[arrow,matrix]{xy}
\usepackage{graphicx,subfigure,tikz,tikz-3dplot}

\usetikzlibrary{positioning,matrix,arrows,calc}
\usepackage{cite}
\usepackage{dsfont}
\usepackage{enumitem}
\usepackage{bbm}
\usepackage{hyperref}

\hypersetup{
	colorlinks=true,
	linkcolor=blue}
\usepackage{array}
\usepackage{xargs}
\usepackage[colorinlistoftodos,prependcaption,textsize=tiny]{todonotes}

\oddsidemargin=0.3in \evensidemargin=0.3in \theoremstyle{plain}
\newtheorem{theorem}{Theorem}[section]
\newtheorem{lemma}[theorem]{Lemma}
\newtheorem{proposition}[theorem]{Proposition}

\theoremstyle{definition}

\newtheorem{definition}[theorem]{Definition}

\numberwithin{equation}{section} 

\newcommand{\R}{\mathbb{R}}
\newcommand{\Prob}{\mathbb{P}}

\newcommand{\C}{\mathbb{C}}
\newcommand{\TV}{\mathrm{TV}}
\newcommand{\Topp}{T_x^{\mathrm{opp}}}
\newcommand{\cU}{\mathcal{U}}
\newcommand{\cQ}{\mathcal{Q}}

\newcommand{\cC}{\mathcal{C}}
\newcommand{\supp}{\operatorname{supp}}

\newcommand{\Leb}{\operatorname{Leb}}

\newcommand{\E}{\mathbb E}
\newcommand{\PP}{\mathcal P}
\newcommand{\TT}{\mathcal T}

\newcounter{hypocounter}
\renewcommand\thehypocounter{(H\arabic{hypocounter})}
\title[Critical Exponents for Random Matrix Products]{\textbf{Critical Exponents and Near-Critical Limits for Random Matrix Products} }
\author{C. Chalhoub}
\address{C. Chalhoub\\ Department of Mathematics,
Imperial College London,
180 Queen's Gate,
London SW7 2AZ,
United Kingdom}
\email{christopher.chalhoub21@imperial.ac.uk}

\author{J.S.W. Lamb}
\address{J.S.W. Lamb\\ Department of Mathematics,
Imperial College London,
180 Queen's Gate,
London SW7 2AZ,
United Kingdom\newline
International Research Center for Neurointelligence (IRCN),
The University of Tokyo, Tokyo 113-0033, Japan}
\email{jsw.lamb@imperial.ac.uk}

\author{M. Rasmussen}
\address{M. Rasmussen\\ Department of Mathematics,
Imperial College London,
180 Queen's Gate,
London SW7 2AZ,
United Kingdom}
\email{m.rasmussen@imperial.ac.uk}

\begin{document}

\begin{abstract}
We study the approach to the zero-Lyapunov-exponent transition for products of independent, identically distributed invertible matrices.
The matrix laws are supported in a common compact subset of \(GL(d,\mathbb R)\) and converge in total variation to a strongly irreducible and proximal law with zero top Lyapunov exponent.
We establish a near-critical renewal theorem for the occupation measure and a functional diffusion limit for the logarithmic cocycle as well as exact asymptotics and scaling limits for three observables describing the transient dynamics. The expected number of escapes
from the unit ball and the expected number of visits to a fixed annulus diverge as the inverse absolute Lyapunov exponent, and both counts have exponential limits after normalization. The time spent on the side opposite to the asymptotic drift has critical exponent two and
converges, after rescaling, to the positive occupation time of Brownian motion with negative drift.
\end{abstract}
\maketitle

\section{Introduction}

Products of random matrices form a classical topic in probability theory, ergodic theory and the theory of random dynamical systems. Since the foundational work of Furstenberg and Kesten
\cite{furstenberg1960products}, they have provided a canonical model for multiplicative stochastic processes. Given an i.i.d.\ sequence $\omega=(A_1,A_2,\ldots)$ of matrices in $GL(d,\mathbb R)$ with common law $\mu$, one studies the products
\[
A_\omega^n=A_n A_{n-1}\cdots A_1.
\]
Under the usual integrability and irreducibility assumptions, for each fixed nonzero initial vector $x$,
\[
\frac1n\log\|A_\omega^n x\|
\longrightarrow\lambda_\mu
\qquad\mathbb P_\mu\text{-almost surely},
\]
where $\lambda_\mu$ is the top Lyapunov exponent. Thus a negative Lyapunov exponent gives exponential convergence to the origin,
whereas a positive exponent gives exponential escape to infinity. The condition $\lambda_\mu=0$ marks the transition between these two regimes.

The Lyapunov exponent, however, describes the eventual exponential trend without quantifying the fluctuations which precede it. In the
contracting regime, an orbit may repeatedly leave a neighbourhood of the origin before contraction prevails. In the expanding regime, it may return many times to a bounded region before escaping
permanently. Near the transition, these transient motions become increasingly pronounced: the logarithmic norm has a small drift, and fluctuations can compete with that drift over long time scales. 
How many times does an orbit cross a fixed threshold before its eventual behaviour takes over? How many visits does it make to a fixed annulus? How much time does it spend on the side opposite to its asymptotic drift?

The asymptotic theory of random matrix products is by now very rich. Furstenberg and Kifer, and Guivarc'h and Raugi, developed the projective and stationary-measure theory under irreducibility and
contraction hypotheses \cite{furstenbergkifer1983,guivarc1985frontiere}; Bougerol and Lacroix gave a systematic account of random matrix products and their applications \cite{RM_Book}. Avila, Eskin and Viana \cite{avila2023continuity} proved continuity of all the Lyapunov exponents when compactly supported matrix laws converge weakly and their supports converge in Hausdorff distance to the support of the limiting law. 
Beyond laws of large numbers, the theory includes central limit theorems, laws of the iterated logarithm, spectral gap estimates and large deviation principles for the logarithmic cocycle; see, for instance, \cite{LePage1982,benoist2016random}. 
More refined local-limit and deviation estimates have been obtained by Xiao, Grama and Liu \cite{XGL_Mod}. The critical regime itself also has an established theory: oscillation at zero Lyapunov exponent was studied by Bougerol \cite{Bougerol1986}, and recurrence and invariant Radon measures for linear random walks were investigated by Guivarc'h and Raja \cite{GuivarchRaja2012}. 
Our question concerns the approach to this recurrent regime from nearby transient systems, and the quantitative laws governing that approach.

This paper develops a theory of critical exponents and near-critical scaling limits for these transient observables. We identify their exact rates of divergence as the Lyapunov exponent approaches zero and determine the limiting distributions after rescaling. 
Two distinct critical scales emerge: crossings and visits to a fixed annulus have exponent $1$, whereas residence on the side opposite to the drift has exponent $2$. The first two observables record what happens near a fixed level; the third also counts time spent far from that level during an excursion against the drift.

To describe the results, fix a compact set $S\subset GL(d,\mathbb R)$ and a strongly irreducible and proximal probability law $\mu_0$ on $S$ with $\lambda_{\mu_0}=0$.
We let $\mu$ approach $\mu_0$ in total variation through probability laws on $S$ with $\lambda_\mu\neq0$. For $x\in\mathbb S^{d-1}$, write
$
S_n(x)=\log\|A_\omega^n x\|.
$
We consider the three observables
\[
\begin{aligned}
N_x
&=\#\{n\geq0:S_n(x)\leq0<S_{n+1}(x)\},\\
L_x(a,b)
&=\#\{n\geq0:a\leq S_n(x)\leq b\},
\qquad a<b,\\
T_x^{\mathrm{opp}}
&=\#\{n\geq0:-\operatorname{sgn}(\lambda_\mu) S_n(x)>0\}.
\end{aligned}
\]
Here $N_x$ counts escapes from the closed unit ball,
$L_x(a,b)$ counts visits to the annulus
$\{y:e^a\leq\|y\|\leq e^b\}$, and
$T_x^{\mathrm{opp}}$ counts the time spent outside the unit ball in the contracting regime and inside the open unit ball in the expanding regime. These observables are almost surely finite for every sufficiently close noncritical measure. At $\mu_0$, the escape count and the number of visits to every fixed annulus are almost surely infinite.

Our main results (Theorems~\ref{MainTheoremEscapesCritical}, \ref{MainTheoremTimeSpentInStripCritical} and~\ref{opposite-exponents}) give the exact asymptotics
\[
\begin{aligned}
\mathbb E_\mu[N_x]
&\sim\frac{\kappa_0}{|\lambda_\mu|},\\
\mathbb E_\mu[L_x(a,b)]
&\sim\frac{b-a}{|\lambda_\mu|},\\
\mathbb E_\mu[T_x^{\mathrm{opp}}]
&\sim\frac{\sigma_0^2}{2|\lambda_\mu|^2}.
\end{aligned}
\]
Here $\sigma_0^2>0$ is the asymptotic variance of the logarithmic cocycle at $\mu_0$, and
\[
\kappa_0
=\frac12
\int_{\mathbf P(\mathbb R^d)}\int_S
|\xi(A,[v])|\,d\mu_0(A)\,d\nu_0([v])>0,
\]
where $\nu_0$ is the $\mu_0$-stationary projective measure and $\xi(A,[v])=\log\|Av\|$ for a unit representative $v$.

The divergence of the means is accompanied by scaling limits. For escapes and annular occupation, we obtain
\[
\frac{|\lambda_\mu|}{\kappa_0}N_x
\xrightarrow{\mathrm d}\operatorname{Exp}(1),
\qquad
\frac{|\lambda_\mu|}{b-a}L_x(a,b)
\xrightarrow{\mathrm d}\operatorname{Exp}(1).
\]
For opposite-phase residence, the limiting law is instead a Brownian occupation time:
\[
|\lambda_\mu|^2T_x^{\mathrm{opp}}
\xrightarrow{\mathrm d}
\mathcal T
:=
\int_0^\infty
\mathbbm{1}_{\{\sigma_0B_t-t>0\}}\,dt,
\]
where $B$ is a standard Brownian motion.

All these limits hold uniformly in the initial direction $x\in\mathbb S^{d-1}$ and allow approach from either sign of the Lyapunov exponent. 
We impose no prescribed one-parameter perturbation: the matrix law may approach $\mu_0$ along any route in the total variation topology within $P(S)$. Thus the powers of $|\lambda_\mu|$ give a common normalization across nearby laws, while the critical law enters the limiting constants through $\kappa_0$ and $\sigma_0^2$. 
To our knowledge, the exact critical asymptotics and scaling limits established here for these random-matrix observables are new.

Several related lines of work place these questions in a broader context. Kesten's work on random difference equations \cite{Kesten1973} initiated the systematic study of heavy tails generated by random matrix products. The associated renewal methods \cite{Kesten1974} have become central to the analysis of affine recursions and their large excursions. 
In particular, Collamore and Mentemeier \cite{Collamore} obtained sharp large-excursion asymptotics and conditioned laws for recursive sequences generated by random matrices. Another related theme is stability near an equilibrium: Arnold \cite{ArnoldStabFormula} connected moment stability and sample stability for linear stochastic systems, while Baxendale and
Khasminskii \cite{baxendale1998stability} studied stability indices and risk probabilities, including probabilities of entering or leaving balls. In the present paper, the spatial thresholds remain fixed while the matrix law approaches the boundary between contraction and expansion. 
The quantities of interest record the repeated crossings and accumulated residence before the eventual drift prevails.

The closest scalar analogue is the small-drift occupation theorem of Schulte-Geers and Stadje \cite{SchulteGeersStadje2017}. For random walks with i.i.d.\ increments of mean $\delta\downarrow0$ and variance tending to $\sigma^2>0$, they show, under a Lindeberg condition, that $\delta^2$ times the number of visits to the negative half-line converges in law to the time that $\sigma B_t+t$ spends there. Taking $\sigma=\sigma_0$ and reflecting the path gives precisely the law $\mathcal T$ above. Small-drift questions for Markov random walks, including renewal estimates and asymptotics for ladder variables, have also been studied in \cite{Fuh2004,Fuh2007}.

For a matrix product, the logarithmic increments need not be independent even though the matrices are: the next increment is $\xi(A_{n+1},[A_\omega^n x])$, whose conditional law generally depends on the direction left by the first $n$ matrices. Changing $\mu$ changes the dynamics of that direction as well as the drift. Thus the scalar theorem does not apply directly, and a limit theorem for one fixed matrix law does not give the uniform estimates needed as $\lambda_\mu\to0$.

There is a further distinction between residence and the two local counts. Under diffusive rescaling a half-line remains a half-line, whereas the logarithmic interval $[a,b]$ shrinks to a point. Convergence of the rescaled paths alone does not determine how many visits the product makes to that interval, or how often a single matrix step takes it across zero. For this we keep track of the projective direction $\theta_n=[A_\omega^n x]$ together with the logarithmic radius, through the renewal measure
\[
U_\mu^x
=
\sum_{n\geq0}
\operatorname{Law}_\mu(\theta_n,S_n(x)).
\]
Theorem~\ref{thm:renewal-limit} gives
\[
|\lambda_\mu|U_\mu^x
\xrightarrow{\mathrm{vague}}
\nu_0\otimes\operatorname{Leb},
\]
uniformly in $x$. Thus the normalized expected occupation becomes stationary in the projective direction and Lebesgue in the
logarithmic radius. Recording the next matrix as well allows us to evaluate the measure of the crossing event and obtain
$\kappa_0$ (Section~\ref{sec:escape-mean}). Evaluating this limit on a fixed strip gives
the coefficient $b-a$ (Section~\ref{sec:strip-visits}). The same theorem also counts visits with a prescribed direction: if $D\subset\mathbf P(\mathbb R^d)$ is Borel, $\nu_0(D)>0$ and $\nu_0(\partial D)=0$, the expected number of visits to the annulus with $\theta_n\in D$ is asymptotic to $(b-a)\nu_0(D)/|\lambda_\mu|$. For $N_x$ and $L_x(a,b)$, a Markov restart argument then turns the uniform asymptotics of the expected remaining counts into the exponential
scaling limits, which yields the proofs of Theorems \ref{MainTheoremEscapesCritical} and \ref{MainTheoremTimeSpentInStripCritical}. Theorem~\ref{thm:renewal-limit} was inspired by the renewal theorem of
Guivarc'h and Le Page \cite[Section~4]{GuivarchLePage2016}, who obtain
a product of a stationary projective measure and Lebesgue measure in
the logarithmic radius, while the matrix law remains fixed in their setting.

The random-walk comparison above explains the origin of the diffusive scale. The accumulated drift has size $n|\lambda_\mu|$, whereas the fluctuations have size $\sqrt n$. 
They become comparable when $n$ is of order $|\lambda_\mu|^{-2}$. We make this balance precise by proving the functional limit (Theorem \ref{lem:diffusive-cocycle-limit})
\[
-\lambda_\mu
S_{\lfloor t/|\lambda_\mu|^2\rfloor}(x)
\xrightarrow{\mathrm d}
\sigma_0B_t-t
\]
on every bounded time interval, in the Skorokhod $J_1$ topology. The sign in the normalization treats the contracting and expanding regimes simultaneously. Occupation of the positive
half-line then corresponds to residence opposite to the original drift. To obtain Theorem~\ref{opposite-exponents}, one must also control the occupation after a large rescaled time, uniformly as the drift vanishes.

The renewal and diffusion limits use the analytic dependence of the transfer operator and its spectral data on the matrix law \cite{AnalyticDependenceIC}. In the renewal argument this operator takes the place of the characteristic function of a single random-walk step: its leading spectral projection near frequency zero gives the stationary factor $\nu_0$ in the limit. In the diffusion argument the quadratic term of the logarithm of its leading eigenvalue gives the variance $\sigma_0^2$; an exponential tilt bounds the expected occupation at late times, allowing passage to the infinite time horizon.

Some further exponents follow directly from the ones above. For example, moving the boundary for opposite-phase residence from the unit sphere to a fixed radius $e^c$, $c\ne0$, changes the count by at most the number of visits to the annulus between the two radii. By Theorem~\ref{MainTheoremTimeSpentInStripCritical} the latter has mean of order $|\lambda_\mu|^{-1}$. The change therefore vanishes in $L^1$ after multiplication by $|\lambda_\mu|^2$: the residence exponent, its leading constant and its limiting law are all unchanged. For escapes, the fluctuations are as large as the mean: the geometric tail bound \eqref{eq:escape-geometric-bounds} makes the squared rescaled counts uniformly integrable, so the exponential limit gives $\operatorname{Var}_\mu(N_x)\sim\kappa_0^2/|\lambda_\mu|^2$.

A different exponent arises if each visit on the opposite side is weighted by its logarithmic distance from the unit sphere. For the sum $\sum_{n\ge0}(-\operatorname{sgn}(\lambda_\mu)S_n(x))_+$, one power of $|\lambda_\mu|$ rescales distance and two rescale time. Multiplying the sum by $|\lambda_\mu|^3$ thus gives the area under the positive part of the rescaled path. The additional point is to control large heights as well as late times: by \eqref{eq:weight-negative-drift}, the probability that the rescaled path exceeds $y\ge0$ at time $t$ is at most $C e^{-c_1t-y/\sigma_0^2}$. This integrable bound lets the same diffusion argument give exponent $3$ for the mean, with leading constant $\mathbb E\int_0^\infty(\sigma_0 B_t-t)_+\,dt$. The renewal and diffusion limits thus provide a starting point for a full critical theory of random matrix products: they distinguish what survives in the Brownian limit from what still depends, as the escape constant $\kappa_0$ does, on the size of a single matrix step.

\section{Setup and Statement of the Main Results}

Let $d\geq 2$. We write $\mathbb{R}^d$ for $d$-dimensional Euclidean space and $\mathbb{S}^{d-1}$ for its unit sphere. For $x\in\mathbb{R}^d$, let
\[
\|x\|=\left(\sum_{i=1}^d |x_i|^2\right)^{1/2}
\]
denote the Euclidean norm. We write $GL(d,\mathbb{R})$ for the group of invertible $d\times d$ real matrices, equipped with the corresponding operator norm
\[
\|A\|=\sup_{x\in\mathbb{S}^{d-1}}\|Ax\|.
\]

Fix a compact set $S\subset GL(d,\mathbb{R})$, and let $P(S)$ denote the set of Borel probability measures on $S$. For $\mu\in P(S)$, set
\[
\Omega=S^{\mathbb{N}},\qquad
\mathcal{F}=\mathcal{B}(S)^{\otimes\mathbb{N}},\qquad
\mathbb{P}_\mu=\mu^{\otimes\mathbb{N}},
\]
and let $\mathbb{E}_\mu$ denote expectation with respect to $\mathbb{P}_\mu$. Write $\omega=(A_1,A_2,\ldots)\in\Omega$, where $A_n$ denotes the $n$th coordinate matrix. Under $\mathbb{P}_\mu$, the matrices $A_1,A_2,\ldots$ are independent with common law $\mu$.

The product at time $n\geq1$ is $A_\omega^n=A_n A_{n-1}\cdots A_1$, with $A_\omega^0=I$. For $x\in \mathbb{R}^d$, $(A_\omega^nx)_{n\geq0}$ is a Markov process on $\mathbb{R}^d$.

We define certain characteristics we impose on the support $W\subset GL(d,\mathbb{R})$ of the measure $\mu$ from \cite{RM_Book}.

\begin{definition}
    The set $W$ is said to be \textit{proximal} if there exists $n\geq1$ and a tuple $(A_1,\dots,A_n)\in W^n$ such that the product $\prod_{i=1}^n A_i$ admits an algebraically simple dominant eigenvalue.
\end{definition}

\begin{definition}
A set $W\subset GL(d,\mathbb{R})$ is called \emph{strongly irreducible}
if there is no finite family of non-zero proper linear subspaces
$V_1,\dots,V_k\subset\mathbb{R}^d$ whose union is preserved by every
matrix in $W$; that is, such that
\[
A\left(\bigcup_{i=1}^k V_i\right)
=
\bigcup_{i=1}^k V_i
\qquad
\text{for every }A\in W.
\]
The weaker condition of \emph{irreducibility} requires only that there
be no non-zero proper linear subspace of $\mathbb{R}^d$ that is
preserved by every matrix in $W$.
\end{definition}

A probability measure $\mu$ on $GL(d,\mathbb{R})$ is called proximal or
strongly irreducible if its support has the corresponding property.

We consider the dynamics of the system on the projective space of directions
\(
\mathbf{P}(\mathbb{R}^d),
\)
where $x\sim y$ if and only if there exists a nonzero $\varsigma\in \R\setminus\{0\}$, such that $x=\varsigma y$. We write $[x]$ for the equivalence class of $x\in \R^d\setminus \{0\}$. For a fixed $x\neq0$, define the induced projective process $(\theta_n)_{n\geq0}$ by
\[
\theta_0=\left[x\right],
\qquad
\theta_n
=
\left[
A_\omega^n x
\right].
\]
Then $(\theta_n)_{n\geq0}$ is a Markov process on
$\mathbf{P}(\mathbb{R}^d)$. By contrast, the radial process
$(\|A_\omega^n x\|)_{n\geq0}$ is generally not Markov on its own,
since its increments depend on the evolving projective direction.

A probability measure $\nu$ on $\mathbf{P}(\mathbb{R}^d)$ is called \emph{$\mu$-stationary} if, for every bounded measurable function
$f:\mathbf{P}(\mathbb{R}^d)\to\mathbb{R}$,
\[
\int_{\mathbf{P}(\mathbb{R}^d)} f([x])\,d\nu([x])
=
\int_S\int_{\mathbf{P}(\mathbb{R}^d)}
f([Ax])\,d\nu([x])\,d\mu(A).
\]

For $A\in S$ and $[x]\in \mathbf{P}(\mathbb{R}^d)$, define \[ \xi(A,[x]) := \log \|A x\|, \]
where $x\in \mathbb{S}^{d-1}$ is a representative. For $x\in\mathbb S^{d-1}$, we can thus write \[ S_n(x) := \log \|A_\omega^n x\|  = \sum_{j=1}^n \xi(A_j, \theta_{j-1}). \]

We next recall the Lyapunov exponents of the random matrix product. For
$M\in GL(d,\mathbb{R})$, let $s_1(M)\geq s_2(M)\geq\cdots\geq s_d(M)>0$ denote its singular values. Since $S$ is compact, the required integrability assumptions are automatic. The Lyapunov exponents associated with $\mu$ are
\[
\lambda_{\mu,i}
:=
\lim_{n\to\infty}
\frac1n
\mathbb{E}_\mu\!\left[\log s_i(A_\omega^n)\right],
\qquad 1\leq i\leq d.
\]
In particular, the top Lyapunov exponent is
\[
\lambda_\mu
:=
\lambda_{\mu,1}
=
\lim_{n\to\infty}
\frac1n
\mathbb{E}_\mu\!\left[\log\|A_\omega^n\|\right].
\]

The Furstenberg--Kesten theorem \cite{furstenberg1960products} gives
\[
\lim_{n\to\infty}
\frac1n\log\|A_\omega^n\|
=
\lambda_\mu
\qquad
\mathbb{P}_\mu\text{-almost surely}.
\]
Furthermore, by Furstenberg and Kifer
\cite{furstenbergkifer1983}, irreducibility implies that for every fixed
$x\in\mathbb{R}^d\setminus\{0\}$,
\[
\lim_{n\to\infty}
\frac1n\log\|A_\omega^n x\|
=
\lambda_\mu
\qquad
\mathbb{P}_\mu\text{-almost surely}.
\]

Strong irreducibility and proximality imply that the top Lyapunov exponent is simple and that the projective process has a unique stationary measure. More precisely, by a theorem of Guivarc'h and Raugi \cite{guivarc1985frontiere}, if $\mu\in P(S)$ is strongly irreducible and proximal, then
\[
\lambda_{\mu,1}>\lambda_{\mu,2},
\]
and there exists a unique $\mu$-stationary probability measure $\nu_\mu$ on $\mathbf{P}(\mathbb{R}^d)$. These properties underlie the spectral estimates used later.

Since $S$ is compact, we may fix $\eta>1$ such that
\[
\max\{\|A\|,\|A^{-1}\|\}\leq\eta
\qquad\text{for every }A\in S.
\]
Consequently,
\[ \eta^{-1}\|x\|\leq\|Ax\|\leq\eta\|x\|\qquad\text{for every }A\in S \text{ and }x\in\mathbb{R}^d.
\]

For $v\in\mathbb{S}^{d-1}$, define the asymptotic variance of the logarithmic cocycle by
\[\sigma_\mu^2:=\lim_{n\to\infty}\frac1n\mathbb{E}_\mu\left[\left(\log\|A_\omega^n v\|-n\lambda_\mu\right)^2\right].\]
Under strong irreducibility, proximality and compact support, this limit is independent of $v$. The central limit theorem for random matrix products, together with its non-degeneracy statement, gives $\sigma_\mu^2>0;$ see \cite[Theorem V.5.1]{RM_Book}, \cite{LePage1982} and \cite[Proposition 3.15(b)]{XGL_Mod}.

We also define the Lyapunov moment function of the random matrix product
(see \cite{ArnoldStabFormula,baxendale1998stability}). For real $q$ in
a neighbourhood of $0$ and $v\in\mathbb{S}^{d-1}$, let
\begin{equation}\label{momlf}
    \Lambda_\mu(q)
    \coloneqq
    \lim_{n\to\infty}
    \frac{1}{n}
    \log
    \mathbb{E}_\mu
    \left[
        e^{q\log\|A_\omega^n v\|}
    \right].
\end{equation}

We are interested in the Lyapunov moment function because it encodes
several important observables of the random matrix product. We recall
some of its basic properties.

\begin{lemma}[\cite{RM_Book}]\label{LyapProps}
Let $S\subset GL(d,\mathbb{R})$ be compact, and let
$\mu\in P(S)$ be strongly irreducible and proximal. Then there exists
$q_0>0$ such that:
\begin{enumerate}
    \item the limit in \eqref{momlf} exists and is independent of
    $v\in\mathbb{S}^{d-1}$ for every real $q$ with $|q|<q_0$;

    \item
    \[
    \Lambda_\mu(0)=0,\qquad
    \Lambda_\mu'(0)=\lambda_\mu,\qquad
    \Lambda_\mu''(0)=\sigma_\mu^2;
    \]

    \item the function $q\mapsto\Lambda_\mu(q)$ is convex and
    real-analytic on $(-q_0,q_0)$;

    \item the function $\Lambda_\mu$ admits a holomorphic extension to
    the complex disc $B_{\mathbb{C}}(0,q_0)$.
\end{enumerate}
\end{lemma}

We use the same notation $\Lambda_\mu$ for this holomorphic extension.
For $k\geq1$, we refer to $\Lambda_\mu^{(k)}(0)$ as the $k$th
Lyapunov moment of the random matrix product.

Let $M(S)$ denote the Banach space of finite complex Borel measures on
$S$, equipped with the total variation norm
\[
\|\rho\|_{\mathrm{TV}}
:=
|\rho|(S)
=
\sup_{\substack{f\in B_b(S)\\ \|f\|_\infty\leq1}}
\left|
\int_S f\,d\rho
\right|,
\]
where $B_b(S)$ denotes the space of bounded Borel measurable
functions on $S$. We equip $P(S)\subset M(S)$ with the metric induced by this norm.

For $\mu\in P(S)$ satisfying the strong irreducibility and proximality
assumptions, we write
$\lambda_\mu,\mathbb{P}_\mu,\mathbb{E}_\mu,\nu_\mu,\Lambda_\mu$ and $\sigma_\mu^2$
respectively for the corresponding top Lyapunov exponent, product measure,
expectation, stationary measure, Lyapunov moment function and
asymptotic variance. We study the behaviour of these
quantities as $\mu\to\mu_0$ in the total variation topology.

For a critical law $\mu_0$, write $\nu_0:=\nu_{\mu_0}$ and
$\sigma_0^2:=\sigma_{\mu_0}^2$. All limits as $\mu\to\mu_0$ below
are taken through probability laws $\mu\in P(S)$ with $\lambda_\mu\ne0$.

\subsection{Scaling Limits and Critical Exponents for Transient Observables}
We now introduce the first observable studied in this paper. For
$x\in\mathbb{S}^{d-1}$ and $\omega\in\Omega$, define the possibly
infinite number of outward crossings of the unit sphere by
\[N_x:=\#\left\{n\geq0:\|A_\omega^n x\|\leq1\text{ and }\|A_\omega^{n+1}x\|>1\right\}.\]

As $\lambda_{\mu}$ depends continuously on $\mu$ in the total variation topology (as stated later in the paper), the condition $\lambda_{\mu_{0}}=0$ implies that $\lambda_{\mu}\to 0$ as $\mu\to\mu_0$. In this setting, we examine the expected number of escapes $N_x$ as $\mu$ approaches the zero-Lyapunov-exponent regime.

\begin{theorem}\label{MainTheoremEscapesCritical}
Let $S\subset GL(d,\mathbb{R})$ be compact, and let
$\mu_0\in P(S)$ be strongly irreducible and proximal with
$\lambda_{\mu_0}=0$. Then
\begin{enumerate}
    \item \begin{equation}\label{Expected_Escapes_Critical_Exponent}
\lim_{\substack{\|\mu-\mu_0\|_{TV}\to 0 \\ \lambda_\mu\ne 0}} \sup_{x\in \mathbb{S}^{d-1}} | |\lambda_\mu| \E_\mu [N_x] - \kappa_0| = 0,
\end{equation}
where
\[ \kappa_0 = \frac{1}{2} \int_{\mathbf{P}(\mathbb{R}^d)} \int_S |\xi(A,[v])| d\mu_0(A) d\nu_0([v]) \in (0,\infty). \]
This establishes an exact critical exponent of $1$ for $\E_\mu [N_x]$ since $\E_\mu [N_x] \sim \kappa_0/|\lambda_\mu|$.
Moreover, $\mathbb{P}_{\mu_0}[N_x=\infty]=1$ for every $x\in\mathbb{S}^{d-1}$.
\item \[ |\lambda_\mu|N_x \xrightarrow{\mathrm{d}} \kappa_0 \mathcal{E} \]
as $\|\mu-\mu_0\|_{TV} \to 0$ uniformly in $x\in \mathbb{S}^{d-1}$, where $\mathcal{E}$ follows an exponential distribution with mean $1$.
\end{enumerate}
\end{theorem}

We can also study the critical behavior of observables describing the time spent in any fixed strip. For any $a<b$, define
\[  L_x(a,b):= \#\{ n\geq 0: \|A_\omega^n x\| \in [e^{a},e^b] \} = \sum_{n\geq 0} \mathbbm{1}_{ \{ a \leq  S_n(x)\leq b \} }. \]

\begin{theorem}\label{MainTheoremTimeSpentInStripCritical}
Let $S\subset GL(d,\mathbb{R})$ be compact, and let
$\mu_0\in P(S)$ be strongly irreducible and proximal with
$\lambda_{\mu_0}=0$. Then for all $a<b$,
\begin{enumerate}
    \item \begin{equation}\label{Expected_Strip_Time_Critical_Exponent}
\lim_{\substack{\|\mu-\mu_0\|_{TV}\to 0 \\ \lambda_\mu\ne 0}} \sup_{x\in \mathbb{S}^{d-1}} | |\lambda_\mu| \E_\mu [L_x(a,b)] - (b-a)| = 0.
\end{equation}
This establishes an exact critical exponent of $1$ for $\E_\mu[L_x(a,b)]$.

Moreover, $\Prob_{\mu_0}[L_x(a,b)=\infty]=1$.
\item \[ |\lambda_\mu|L_x(a,b) \xrightarrow{\mathrm{d}} (b-a)  \mathcal{E} \]
as $\|\mu-\mu_0\|_{TV} \to 0$ uniformly in $x\in \mathbb{S}^{d-1}$, where $\mathcal{E}$ follows an exponential distribution with mean $1$.
\end{enumerate}

\end{theorem}

We finally present an observable with critical exponent $2$. For $x\in \mathbb{S}^{d-1}$ and $\lambda_\mu\ne 0$, define the opposite-phase residence \[ \Topp := \#\{ n\geq 0: -\operatorname{sgn}(\lambda_\mu) S_n(x)>0 \}. \]
Thus $\Topp$ counts the time spent outside the unit ball when $\lambda_\mu<0$, and inside the open unit ball when $\lambda_\mu>0$. The boundary is excluded in both cases. The drift $n\lambda_\mu$ and the fluctuations of order $\sqrt n$ balance when $n$ is of order $|\lambda_\mu|^{-2}$. On this time scale, the logarithmic norm, multiplied by $-\lambda_\mu$, approaches Brownian motion with drift $-1$.

\begin{theorem}\label{opposite-exponents}
Let $S\subset GL(d,\mathbb{R})$ be compact, and let
$\mu_0\in P(S)$ be strongly irreducible and proximal with
$\lambda_{\mu_0}=0$. Then,
\begin{enumerate}
    \item \begin{equation}\label{eq:opposite-residence-limit}
\lim_{\substack{\|\mu-\mu_0\|_{TV}\to 0 \\ \lambda_\mu\ne 0}} \sup_{x\in\mathbb S^{d-1}}
\left|\,|\lambda_\mu|^2\E_\mu[\Topp]-\frac{\sigma_0^2}{2}\right|
=0.
\end{equation}
This establishes an exact critical exponent of $2$ for $\E_\mu[\Topp]$.
\item \begin{equation}\label{eq:opposite-residence-law}
|\lambda_\mu|^2\Topp\ \xrightarrow{\ d\ }\
\TT:=\int_0^\infty\mathbbm{1}_{\{\sigma_0 B_t-t>0\}}\,dt,
\end{equation}
as $\|\mu-\mu_0\|_{TV} \to 0$ uniformly in $x\in \mathbb{S}^{d-1}$, where $B$ is a standard Brownian motion.
\end{enumerate}
\end{theorem}

\subsection{Near-critical renewal and diffusion limits}

The following two results describe the near-critical process itself
and underlie the preceding observable limits.
For $\lambda_\mu\ne0$ and $x\in\mathbb S^{d-1}$, define
the renewal measure by
\begin{equation}\label{eq:renewal-measure}
U_\mu^x(B)
:=\sum_{n=0}^{\infty}
\Prob_\mu\bigl((\theta_n,S_n(x))\in B\bigr),
\qquad B\in\mathcal B(\mathbf P(\R^d)\times\R).
\end{equation}
This measure counts expected visits to each set.

\begin{theorem}[Near-critical renewal limit]\label{thm:renewal-limit}
Let $S\subset GL(d,\mathbb{R})$ be compact, and let
$\mu_0\in P(S)$ be strongly irreducible and proximal with
$\lambda_{\mu_0}=0$.
For all noncritical $\mu\in P(S)$ sufficiently close to $\mu_0$,
the measures $U_\mu^x$ are locally finite. Moreover, for every $F\in C_c(\mathbf P(\R^d)\times\R)$,
\begin{equation}\label{eq:renewal-limit}
\lim_{\substack{\mu\to\mu_0\\\lambda_\mu\ne0}}
\sup_{x\in\mathbb S^{d-1}}
\left|
|\lambda_\mu|\int F\,dU_\mu^x
-\int_{\mathbf P(\R^d)}\int_\R F([v],s)\,ds\,d\nu_0([v])
\right|=0.
\end{equation}
Thus $|\lambda_\mu|U_\mu^x\to\nu_0\otimes\Leb$ vaguely,
uniformly in $x\in\mathbb S^{d-1}$.
\end{theorem}

On the diffusive time scale, define the process
\begin{equation}\label{eq:diffusive-cocycle-process}
Z_\mu^x(t):=-\lambda_\mu
S_{\lfloor t/|\lambda_\mu|^2\rfloor}(x),\qquad t\ge0.
\end{equation}
For $K>0$, let $D([0,K])$ be the Skorokhod space of c\`adl\`ag functions on $[0,K]$ with values in $\R$.

\begin{theorem}[Diffusive cocycle limit]\label{lem:diffusive-cocycle-limit}
Let $S\subset GL(d,\mathbb{R})$ be compact, and let
$\mu_0\in P(S)$ be strongly irreducible and proximal with
$\lambda_{\mu_0}=0$.
For every $0<K<\infty$, along every sequence of noncritical laws
$\mu_j\to\mu_0$ and directions $x_j\in\mathbb S^{d-1}$,
\[
Z_{\mu_j}^{x_j}\ \xrightarrow{\ d\ }\
(\sigma_0B_t-t)_{0\le t\le K}
\quad\text{in }D([0,K])\text{ with the }J_1\text{ topology},
\]
where $B$ is a standard Brownian motion on $\R$.
\end{theorem}

\section{Organization of the Paper}\label{organization}
Section~\ref{PerturbationTheorySection} establishes the uniform
spectral estimates used in
Section~\ref{Summed_Transfer_Operator} to analyse the summed transfer
operator. Section~\ref{renewal} proves the near-critical renewal limit
of Theorem~\ref{thm:renewal-limit}. Sections~\ref{sec:escape-mean}
and~\ref{sec:strip-visits} use this limit and a Markov restart argument
to prove Theorems~\ref{MainTheoremEscapesCritical}
and~\ref{MainTheoremTimeSpentInStripCritical}, respectively.
Finally, Section~\ref{sec:residence} proves the diffusive cocycle limit
of Theorem~\ref{lem:diffusive-cocycle-limit} and applies it to obtain
Theorem~\ref{opposite-exponents}.

\section{Perturbation Theory}\label{PerturbationTheorySection}

Throughout the following sections, fix a strongly irreducible and proximal
law $\mu_0\in P(S)$ with $\lambda_{\mu_0}=0$.

For $\alpha\in(0,1)$, let
$\mathcal{C}^\alpha(\mathbf{P}(\mathbb{R}^d))$
denote the Banach space of complex-valued $\alpha$-Hölder continuous
functions on $\mathbf{P}(\mathbb{R}^d)$. We equip
$\mathbf{P}(\mathbb{R}^d)$ with the metric
\[
d_{\mathbf{P}}([x],[y])
:=
\min\{\|x-y\|,\|x+y\|\},
\qquad x,y\in\mathbb{S}^{d-1},
\]
and endow $\mathcal{C}^\alpha(\mathbf{P}(\mathbb{R}^d))$ with the norm
\[
\|f\|_{\mathcal C^\alpha}
=
\|f\|_\infty+[f]_\alpha,
\qquad
[f]_\alpha
=
\sup_{[x]\neq[y]}
\frac{|f([x])-f([y])|}
     {d_{\mathbf{P}}([x],[y])^\alpha}.
\]
Let
$\mathbbm{1}_{\mathbf{P}(\mathbb{R}^d)}$
denote the constant function equal to $1$ on
$\mathbf{P}(\mathbb{R}^d)$. We write
\(
L\bigl(\mathcal{C}^\alpha(\mathbf{P}(\mathbb{R}^d))\bigr)
\)
for the space of bounded linear operators on $\mathcal{C}^\alpha(\mathbf{P}(\mathbb{R}^d))$ and
\(
\mathcal{C}^\alpha(\mathbf{P}(\mathbb{R}^d))^*
\)
for its dual space. We use $\|\cdot\|_{\mathcal C^\alpha}$ also for
the induced operator norm.

For $q\in\mathbb{C}$ and $\mu\in M(S)$, define the
\emph{$q$-twisted Koopman operator}
\(
\mathcal{P}_{q,\mu}
\in
L\bigl(\mathcal{C}^\alpha(\mathbf{P}(\mathbb{R}^d))\bigr)
\)
by
\[
\mathcal{P}_{q,\mu}f([v])
=
\int_S
e^{q\log\|Av\|}
f([Av])\,d\mu(A),
\qquad
f\in\mathcal{C}^\alpha(\mathbf{P}(\mathbb{R}^d)),
\]
where $v\in\mathbb{S}^{d-1}$ is any unit representative of $[v]$.
This definition is independent of the choice of unit representative.

For $\mu\in P(S)$, we have \[ \mathcal{P}_{q,\mu}^n f ([x]) = \mathbb{E}_{\mu} \left[ e^{qS_n(x)} f(\theta_n)  \right] .\]

The following standard spectral result connects the Lyapunov moment
function with the twisted Koopman operator.

\begin{proposition}[{\cite[V.4.3]{RM_Book}}]
\label{Spectral_Gap_Operator}
Let $S\subset GL(d,\mathbb{R})$ be compact, and let
$\mu\in P(S)$ be strongly irreducible and proximal. There exists
$\alpha_0>0$ such that, for every $\alpha\in(0,\alpha_0)$, there exists
$q_0>0$ such that, for every $q\in B_{\mathbb{C}}(0,q_0)$, the operator
\[
\mathcal{P}_{q,\mu}
:
\mathcal{C}^\alpha(\mathbf{P}(\mathbb{R}^d))
\longrightarrow
\mathcal{C}^\alpha(\mathbf{P}(\mathbb{R}^d))
\]
has an isolated simple eigenvalue $e^{\Lambda_\mu(q)}$.
\end{proposition}

Henceforth, fix $\alpha>0$ sufficiently small for the spectral gap
and the nonzero-frequency spectral estimate used below at $\mu_0$.

We recall the definition of Fréchet-analyticity (holomorphicity) for maps between Banach spaces following \cite{chae2020holomorphy}. It generalizes the usual definition on $\mathbb{C}$, characterized by the existence of a power series expansion on $\mathbb{C}$.
In the following, let $E$ and $F$ be complex Banach spaces and $U\subset E$ an open set.

\begin{definition}
     The map $T: U\to F$ is said to be \textit{Fréchet-holomorphic} (or Fréchet-analytic) at $x_0\in U$ if there exist $r>0$ with $B_E(x_0,r)\subset U$ and a sequence of continuous, symmetric complex $n$-linear maps $T_n: E^n\to F$ for $n \geq 0$ such that
    \[ T(x) = \sum_{n=0}^\infty T_n (x-x_0,\dots,x-x_0) \]
    for all $x\in B(x_0,r)$, where the series above converges in the norm of $F$. $T$ is said to be Fréchet-holomorphic on $U$ if it is Fréchet-holomorphic at every $x_0\in U$.
\end{definition}

We typically refer to Fréchet-holomorphicity when simply writing holomorphic or analytic in this paper.
\begin{proposition}[Uniform analytic spectral data]\label{prop:analytic-spectral-data}
There exist $r_0>0$, a total-variation neighbourhood $\cU_0\subset M(S)$ of $\mu_0$, constants $C_0<\infty$ and $\rho_0\in(0,1)$, and analytic maps
\[(q,\mu)\longmapsto\Lambda_\mu(q) \in \C,\qquad (q,\mu)\longmapsto\Pi_{q,\mu}\in L\bigl(\mathcal{C}^\alpha(\mathbf{P}(\R^d))\bigr),\]
defined for $|q|<r_0$ and $\mu\in\cU_0$, such that, with $\mathcal Q_{q,\mu}:=\mathcal P_{q,\mu}(I-\Pi_{q,\mu})$,
\begin{equation}\label{eq:uniform-spectral-decomp}
\mathcal P_{q,\mu}^{\,n}=e^{n\Lambda_\mu(q)}\Pi_{q,\mu}+\mathcal Q_{q,\mu}^{\,n},
\qquad
\|\mathcal Q_{q,\mu}^{\,n}\|_{\mathcal C^\alpha}\leq C_0\rho_0^n,
\end{equation}
for every $n\geq1$. For every $\mu\in P(S)\cap\cU_0$,
\begin{equation}\label{eq:projection-at-zero}
\Pi_{0,\mu}f=\nu_\mu(f) \mathbbm{1}_{\mathbf{P}(\mathbb{R}^d)} ,
\end{equation}
where $\nu_\mu$ is the unique $\mu$-stationary probability measure on $\mathbf{P}(\R^d)$. Moreover, $(q,\mu)\mapsto\Pi_{q,\mu}$ and $q\mapsto\Lambda_\mu(q)$ have uniformly bounded $q$-derivatives of every fixed finite order on smaller product neighbourhoods.
\end{proposition}

\begin{proof}
We identify the projective Hölder space with the closed invariant subspace of even functions in the Hölder space on the sphere.\par
Under this identification, the operator-analyticity argument of \cite[Proposition 4.1]{AnalyticDependenceIC} applies to our projective operator. The results of \cite{AnalyticDependenceIC} yield joint analyticity in $q$ and $\mu$ in an open neighbourhood $ B_{r_0}(0) \times \cU_0 $ of $(0,\mu_0)$ of the maps $\Lambda_\mu(q)$ and $\mathcal P_{q,\mu}$. Moreover, defining \[ \Pi_{q,\mu} := \frac{1}{2\pi i} \int_{\Gamma} (\xi I - \mathcal{P}_{q,\mu})^{-1} d\xi \]
as the Riesz projection of the transfer operator, where $\Gamma$ is a fixed
positively oriented contour enclosing $1$ and no other spectral point of
$\mathcal P_{0,\mu_0}$ and remains isolating after shrinking the neighbourhood,
the argument of \cite[Proposition 5.1]{AnalyticDependenceIC} similarly yields
joint analyticity of $\Pi_{q,\mu}$ in $q$ and $\mu$ in the neighbourhood above.

Now set \[\mathcal{Q}_{q,\mu} = \mathcal{P}_{q,\mu} (I - \Pi_{q,\mu}).\] The definition implies that the map $(q,\mu)\mapsto \mathcal{Q}_{q,\mu}$ is also analytic.
By \cite[Theorem V.4.3]{RM_Book}, there exists $\rho_0>0$ such that $\text{specrad}(\mathcal{Q}_{0,\mu_0}) < \rho_0 <1$ (where specrad refers to the spectral radius of the operator). The defining properties of the Riesz projection give
\[
\mathcal P_{q,\mu}\Pi_{q,\mu}=e^{\Lambda_\mu(q)}\Pi_{q,\mu},
\qquad
\Pi_{q,\mu}\mathcal Q_{q,\mu}=\mathcal Q_{q,\mu}\Pi_{q,\mu}=0,
\]
and the decomposition
\[ \mathcal{P}_{q,\mu} = e^{\Lambda_\mu (q)} \Pi_{q,\mu} + \mathcal{Q}_{q,\mu}. \]
It follows that
\[ \mathcal{P}_{q,\mu} ^n  = e^{n\Lambda_\mu (q)} \Pi_{q,\mu} + \mathcal{Q}^n_{q,\mu}, \]
and Dunford's formula gives the representation
\[ \mathcal{Q}^n_{q,\mu} = \frac{1}{2\pi i} \int_{|\xi|=\rho_0} \xi^n (\xi I - \mathcal{Q}_{q,\mu})^{-1} d\xi   \]
for all $(q,\mu)$ in the neighbourhood $B_{r_0}(0)\times \cU_0$ (after potentially shrinking it). By continuity of $\mathcal{Q}_{q,\mu}$, there exists $M>0$ such that
\[ \sup_{|q|<r_0, \mu \in \cU_0} \sup_{|z|=\rho_0} \|(zI - \mathcal{Q}_{q,\mu})^{-1} \|_{\mathcal C^\alpha} \leq M. \]
This implies that
\[ \|\mathcal{Q}^n_{q,\mu}\|_{\mathcal C^\alpha} \leq M \rho_0^{n+1}, \]
where $M$ is independent of $q,\mu$ and $n$, as desired.

Since strong irreducibility and proximality jointly form an open
condition in the weak topology, and hence in the total variation
topology \cite[Section~1]{GuivarchLePage2016}, we may shrink $\cU_0$
so that every $\mu\in P(S)\cap\cU_0$ is strongly irreducible and proximal, and \cite[Theorem V.4.3(i)]{RM_Book} yields $\Pi_{0,\mu}f=\nu_\mu(f) \mathbbm{1}_{\mathbf{P}(\mathbb{R}^d)}$.
Analyticity also yields the following bounded q-derivatives (which we require for Taylor expansions): for all $k\in \mathbb{N}$, there exist open neighbourhoods $\cU_k$ of $\mu_0$ and $r_k>0$ such that
\[ \sup_{\mu \in \cU_k, |q|<r_k} ( \|\partial^k_q \Pi_{q,\mu}\|_{\mathcal C^\alpha} + |\partial^k_q \Lambda_\mu(q)| )<\infty, \]
where the operator derivatives correspond to Fréchet-derivatives.

\end{proof}

Since $\lambda_\mu = \Lambda'_\mu(0)$ and $\sigma^2_\mu = \Lambda_\mu''(0)$ we have $\lambda_\mu \to 0$ and $\sigma^2_\mu \to \sigma^2_{\mu_0}>0$ as $\|\mu-\mu_0\|_{TV} \to 0$.

\section{The Summed Transfer Operator}\label{Summed_Transfer_Operator}

We will study expected occupation times by summing over all times before
taking the limit $\mu\to\mu_0$. The transfer-operator identity is
\begin{equation}\label{eq:operator-Pn}
\E_\mu\!\left[e^{zS_n(x)}f(\theta_n)\right]
=\bigl(\PP_{z,\mu}^{\,n}f\bigr)([x]),
\qquad f\in\mathcal C^\alpha(\mathbf P(\R^d)),
\end{equation}
where $x$ is a unit vector, $S_n(x)=\log\|A_\omega^nx\|$ and
$\theta_n=[A_\omega^nx]$. Thus, whenever the operator series converges,
\begin{equation}\label{eq:summed-identity}
\sum_{n=0}^\infty
\E_\mu\!\left[e^{zS_n(x)}f(\theta_n)\right]
=\left(\sum_{n=0}^\infty\PP_{z,\mu}^{\,n}f\right)([x])
=\bigl((I-\PP_{z,\mu})^{-1}f\bigr)([x]).
\end{equation}
For $z=q+it$, the real part $q$ will be associated with a function of $\lambda_\mu$ and $t$ is
the Fourier frequency. The weight $q$ will make the time sum finite, and the
inverse in \eqref{eq:summed-identity} then lets us analyse the whole sum
at once. We separate frequencies into frequencies away from zero where the resolvent operator remains bounded, and those near zero  where we observe a critical divergence.

\begin{lemma}[Frequencies away from zero]\label{lem:away-frequency}
Fix $0<a<T<\infty$. There exist a neighbourhood $\cU_{a,T}$ of
$\mu_0$, numbers $h>0$ and $\rho\in(0,1)$, and $C<\infty$ such that
\begin{equation}\label{eq:away-frequency}
\|\PP_{b+it,\mu}^{\,n}\|_{\mathcal C^\alpha}\le C\rho^n
\quad\text{if }\mu\in P(S)\cap\cU_{a,T},\ |b|\le h,
\ a\le|t|\le T,\ n\ge0.
\end{equation}
Consequently,
\begin{equation}\label{eq:away-resolvent}
\sum_{n=0}^\infty\PP_{b+it,\mu}^{\,n}
=(I-\PP_{b+it,\mu})^{-1},\qquad
\|(I-\PP_{b+it,\mu})^{-1}\|_{\mathcal C^\alpha}\le C',
\end{equation}
with convergence in operator norm.
\end{lemma}

\begin{proof}

Under our compact-support, strong irreducibility and proximality assumptions, \cite[Proposition 3.7]{XGL_Mod} applied at $s=0$ gives
\begin{equation}\label{eq: specPit}
    \text{specrad}(\mathcal{P}_{it,\mu_0}) <1 \quad \text{for all } t\neq 0.
\end{equation}

Set $K_{a,T}:=[-T,-a]\cup[a,T]$. By upper semicontinuity of the spectral radius and compactness of
the reference family, operator continuity allows us to choose
$0<h\le1$ and a neighbourhood $\cU_{a,T}$ of $\mu_0$ such that
\[
\operatorname{spec}(\PP_{b+it,\mu})
\subset \{w\in\mathbb C:|w|<\rho\}
\]
for some $\rho<1$ whenever $|b|\le h$, $t\in K_{a,T}$ and
$\mu\in P(S)\cap\cU_{a,T}$.

We next bound the resolvents uniformly on the circle $|w|=\rho$.
Since inversion is continuous on the set of invertible operators,
the map
\[
(w,t)\longmapsto R_t(w):=(wI-\PP_{it,\mu_0})^{-1}
\]
is continuous on the compact set
$\{|w|=\rho\}\times K_{a,T}$. Hence
\[
M_0:=\sup_{\substack{|w|=\rho\\t\in K_{a,T}}}
\|R_t(w)\|_{\mathcal C^\alpha}<\infty.
\]
Set $E_{b,t,\mu}:=\PP_{b+it,\mu}-\PP_{it,\mu_0}$.
Operator continuity gives us
\begin{equation}\label{eq:operator-continuity}
\|\PP_{z,\mu}-\PP_{w,\nu}\|_{\mathcal C^\alpha}
\le C_T\bigl(|z-w|+\|\mu-\nu\|_{\TV}\bigr),
\qquad |z|,|w|\le T,\quad \mu,\nu\in P(S).
\end{equation}
Applying \eqref{eq:operator-continuity} on the disc of radius
$T+1$, and shrinking $h$ and $\cU_{a,T}$ if necessary, gives
\[
\|E_{b,t,\mu}\|_{\mathcal C^\alpha}
\le C_{T+1}\bigl(|b|+\|\mu-\mu_0\|_{\TV}\bigr)
\le \frac{1}{2M_0}.
\]
The factorization
\[
wI-\PP_{b+it,\mu}
=\bigl(I-E_{b,t,\mu}R_t(w)\bigr)
  (wI-\PP_{it,\mu_0})
\]
and the Neumann series therefore yield
\[
(wI-\PP_{b+it,\mu})^{-1}
=R_t(w)\sum_{k\ge0}\bigl(E_{b,t,\mu}R_t(w)\bigr)^k,
\qquad
\|(wI-\PP_{b+it,\mu})^{-1}\|_{\mathcal C^\alpha}\le 2M_0.
\]
Thus, writing $M=2M_0$, Dunford's formula (see for example \cite[Section V.3, p.~112]{RM_Book}) on the circle $|w|=\rho$ gives
\[
\PP_{b+it,\mu}^{\,n}
=\frac{1}{2\pi i}\int_{|w|=\rho}
w^n(wI-\PP_{b+it,\mu})^{-1}\,dw,
\]
and consequently
\[
\|\PP_{b+it,\mu}^{\,n}\|_{\mathcal C^\alpha}
\le \frac{1}{2\pi}(2\pi\rho)\rho^n M
=M\rho^{n+1}.
\]
Since $\rho<1$, the series of powers converges in operator norm
and satisfies
\[
(I-\PP_{b+it,\mu})^{-1}
=\sum_{n\ge0}\PP_{b+it,\mu}^{\,n},
\qquad
\|(I-\PP_{b+it,\mu})^{-1}\|_{\mathcal C^\alpha}
\le \frac{M\rho}{1-\rho}.
\]
This proves \eqref{eq:away-resolvent}.

\end{proof}

We now move to the critical scaling. The second-order approximation to $\Lambda_\mu(q)$ is
$\lambda_\mu q+\tfrac12\sigma_0^2q^2$. Its minimum occurs at
\begin{equation}\label{eq:exponential-weight}
q_\mu:=-\frac{\lambda_\mu}{\sigma_0^2},
\end{equation}
where its value is $-\lambda_\mu^2/(2\sigma_0^2)$. This motivates the
weight: for nearby noncritical laws, $\Lambda_\mu(q_\mu)<0$, so the
weighted expectations decay exponentially in time. The same quadratic
approximation will determine the critical limit.

\begin{lemma}[The critical resolvent]\label{lem:critical-resolvent}
There exist a neighbourhood $\cU\subset\cU_0$ of $\mu_0$ and constants
$r,c,C>0$ such that the following hold for
$\mu\in P(S)\cap\cU$ with $\lambda_\mu\ne0$.

\begin{enumerate}[label=\textnormal{(\roman*)},leftmargin=*,itemsep=6pt]
\item Uniformly in the initial direction $x$,
\begin{equation}\label{eq:weight-negative-drift}
\E_\mu[e^{q_\mu S_n(x)}]
\le Ce^{-c|\lambda_\mu|^2n},
\qquad n\ge0.
\end{equation}

\item For $|t|\le r$, the operator series converges in norm and
\begin{equation}\label{eq:resolvent-decomposition}
\sum_{n=0}^\infty\PP_{q_\mu+it,\mu}^{\,n}
=(I-\PP_{q_\mu+it,\mu})^{-1}.
\end{equation}
Moreover,
\begin{equation}\label{eq:resolvent-denominator}
\|(I-\PP_{q_\mu+it,\mu})^{-1}\|_{\mathcal C^\alpha}
\le\frac{C}{|\lambda_\mu|^2+t^2}.
\end{equation}

\item For every $M>0$, as $\mu\to\mu_0$ through laws with
$\lambda_\mu\ne0$,
\begin{equation}\label{eq:critical-resolvent-limit}
\sup_{|u|\le M}
\left\|
|\lambda_\mu|^2(I-\PP_{q_\mu+i|\lambda_\mu|u,\mu})^{-1}
-\frac{2\sigma_0^2}{1+\sigma_0^4u^2}\Pi_{0,\mu_0}
\right\|_{\mathcal C^\alpha}
\longrightarrow0.
\end{equation}
\end{enumerate}
\end{lemma}

\begin{proof}
In a neighbourhood of $z=0$ we can write \begin{equation}\label{eq:complex-Lambda-Taylor}
\Lambda_\mu(z)=\lambda_\mu z+\frac12\sigma_\mu^2z^2+O(|z|^3),
\end{equation}
where the remainder is uniform in $\mu$ (from the uniform q-derivatives results of Proposition \ref{prop:analytic-spectral-data}). Setting $z=q_\mu + it$,
the quadratic expression using the $\mu_0$-variance
simplifies exactly:
\[
\lambda_\mu(q_\mu+it)
+\frac12\sigma_0^2(q_\mu+it)^2
=-\frac{\lambda_\mu^2}{2\sigma_0^2}
-\frac12\sigma_0^2t^2.
\]
The imaginary cross terms cancel because
$\lambda_\mu+\sigma_0^2q_\mu=0$. Hence \eqref{eq:complex-Lambda-Taylor}
becomes
\begin{equation}\label{eq:completed-square}
\begin{split}
\Lambda_\mu(q_\mu+it)
={}&-\frac{\lambda_\mu^2}{2\sigma_0^2}-\frac12\sigma_0^2t^2\\
&+O\!\left(
\bigl(|\sigma_\mu^2-\sigma_0^2|+|\lambda_\mu|+|t|\bigr)
\bigl(|\lambda_\mu|^2+t^2\bigr)\right).
\end{split}
\end{equation}
The error is uniform. Since $\sigma_\mu^2\to\sigma_0^2$ and
$\lambda_\mu\to0$, choosing $r$ and $\cU$ sufficiently small yields
\begin{equation}\label{eq:vertical-quadratic-bound}
\Re\Lambda_\mu(q_\mu+it)
\le-c\bigl(|\lambda_\mu|^2+t^2\bigr),\qquad |t|\le r.
\end{equation}

The proof of (i) simply uses \eqref{eq:operator-Pn} and the spectral decomposition established in the previous section: for all $n\geq 1,$
\begin{align*}
\E_\mu[e^{q_\mu S_n(x)}]
&=e^{n\Lambda_\mu(q_\mu)}
  (\Pi_{q_\mu,\mu} \mathbbm{1}_{\mathbf{P}(\mathbb{R}^d)} )([x])
  +(\cQ_{q_\mu,\mu}^{\,n} \mathbbm{1}_{\mathbf{P}(\mathbb{R}^d)} )([x])\\
&\le Ce^{-c|\lambda_\mu|^2n}+C\rho_0^n
\le C'e^{-c'|\lambda_\mu|^2n}.
\end{align*}
for some uniform constants $C'$ and $c'$.
Equation \eqref{eq:vertical-quadratic-bound} gives us, for $z=q_\mu + it$,
\[\left|e^{\Lambda_\mu(z)}\right| = e^{\Re \Lambda_\mu(z)} <1.\]
Hence, both terms in the spectral decomposition are summable in operator norm:
\begin{align*}
(I-\PP_{z,\mu})^{-1}=\sum_{n=0}^\infty\PP_{z,\mu}^{\,n}
&=I+\sum_{n=1}^\infty
\left(e^{n\Lambda_\mu(z)}\Pi_{z,\mu}+\cQ_{z,\mu}^{\,n}\right)\\
&=\frac{\Pi_{z,\mu}}{1-e^{\Lambda_\mu(z)}}
+\underbrace{(I-\Pi_{z,\mu})+
\sum_{n=1}^\infty\cQ_{z,\mu}^{\,n}}_{\displaystyle\mathcal R_{z,\mu}},
\end{align*}
and the uniform bounds on
$\Pi_{z,\mu}$ and $\cQ_{z,\mu}^{\,n}$ (Proposition \ref{prop:analytic-spectral-data}) give
$\|\mathcal R_{z,\mu}\|_{\mathcal C^\alpha}\le C$.
Furthermore,
\[
|1-e^{\Lambda_\mu(z)}|
\ge1-e^{\Re\Lambda_\mu(z)}
\ge1-e^{-c(|\lambda_\mu|^2+t^2)}
\ge c'(|\lambda_\mu|^2+t^2).
\]
This yields the proof of (ii) via a triangle inequality:
\[ \|(I- \PP_{z,\mu})^{-1} \|_{\mathcal C^\alpha} \leq \frac{\|\Pi_{z,\mu}\|_{\mathcal C^\alpha}}{|1-e^{\Lambda_\mu(z)}|} + \|\mathcal R_{z,\mu}\|_{\mathcal C^\alpha} \leq \frac{C}{|\lambda_\mu|^2+t^2} . \]

Fix $M>0$. We now prove the critical limit (iii) by resolving the critical frequency scale: set $t=u|\lambda_\mu|$. Substituting this in \eqref{eq:completed-square} gives
\[
\frac{\Lambda_\mu(q_\mu+i|\lambda_\mu|u)}{|\lambda_\mu|^2}
=-\frac{1}{2\sigma_0^2}-\frac12\sigma_0^2u^2+o(1),
\]
where the error tends to zero uniformly in $u\in [-M,M]$. Then, using $1-e^w=-w+O(|w|^2)$, we obtain
\begin{equation}\label{eq:scaled-scalar-resolvent}
\frac{|\lambda_\mu|^2}
{1-e^{\Lambda_\mu(q_\mu+i|\lambda_\mu|u)}}
\longrightarrow\frac{2\sigma_0^2}{1+\sigma_0^4u^2}
\end{equation}
uniformly on $[-M,M]$ as $\|\mu - \mu_0\|_{TV}\to 0$.
Using
\[
|\lambda_\mu|^2(I-\PP_{z,\mu})^{-1}
=
\frac{|\lambda_\mu|^2}{1-e^{\Lambda_\mu(z)}}\Pi_{z,\mu}
+|\lambda_\mu|^2\mathcal R_{z,\mu},
\]
since the Riesz projection converges $\Pi_{q_\mu+i|\lambda_\mu|u,\mu} \to \Pi_{0,\mu_0}$ in operator norm, this proves (iii).
\end{proof}

\section{Near-Critical Renewal Measure}\label{renewal}

Throughout this section, $\mu\in P(S)\cap\cU_0$, with $\cU_0$ reduced
so that Lemma~\ref{lem:critical-resolvent} applies.
For $\lambda_\mu\ne0$ and $x\in\mathbb S^{d-1}$, recall
the renewal measure $U_\mu^x$ defined in
\eqref{eq:renewal-measure}. Our aim is to prove
Theorem~\ref{thm:renewal-limit}, which states that
$|\lambda_\mu|U_\mu^{x}$ converges locally to $\nu_0\otimes\Leb$.
The method is to Fourier-transform the logarithmic coordinate and sum
the resulting transfer operators using Section \ref{Summed_Transfer_Operator}.

The time sum in \eqref{eq:renewal-measure} has infinite total mass.
To make it summable, use $q_\mu=-\lambda_\mu/\sigma_0^2$ from
\eqref{eq:exponential-weight} and set
\begin{equation}\label{eq:weighted-renewal-measure}
W_\mu^x(d[v],ds)
:=|\lambda_\mu|e^{q_\mu s}U_\mu^x(d[v],ds).
\end{equation}
The expectation bound \eqref{eq:weight-negative-drift} gives
\begin{equation}\label{eq:weighted-total-mass}
W_\mu^x(\mathbf P(\R^d)\times\R)
=|\lambda_\mu|
\sum_{n\ge0}\E_\mu[e^{q_\mu S_n(x)}]
\le \frac{C}{|\lambda_\mu|}<\infty.
\end{equation}
Thus $W_\mu^x$ is finite and $U_\mu^x$ is locally finite;
both are Radon measures. Since $q_\mu\to0$, the weight tends to one
uniformly on compact $s$-intervals. We will first prove the limit for
$W_\mu^x$ and then remove the weight.

Renewal measures for Markov additive processes and random matrix
products appear in \cite{Kesten1974,Kesten1973,LePage1983}, and related
uniform and small-drift results are developed in
\cite{Fuh2004,Fuh2007,BlanchetGlynn2007}. The stationary-measure/Lebesgue
product also occurs in the fixed-law renewal theorems of
\cite[Section 4]{GuivarchLePage2016}, where the starting radius tends
to zero. The difference here is that the law approaches a critical law while the starting
logarithmic radius remains bounded.

Recall that vague convergence means convergence against every
continuous, compactly supported test function. The following consequence
will later allow us to use crossing indicators. Let $\operatorname{Disc}(G)$ denote the set of discontinuities of a function $G$.

\begin{lemma}[Portmanteau consequence]\label{lem:vague-portmanteau}
Let $(m_j)_{j\geq 1}$ and $m$ be Radon measures such that $m_j\to m$ vaguely on a locally compact metric
space. If $G$ is bounded, Borel measurable, vanishes outside a compact
set, and $m(\operatorname{Disc}(G))=0$, then
$\int G\,dm_j\to\int G\,dm$. In particular, $m_j(B)\to m(B)$ if
$B$ is contained in a compact set and $m(\partial B)=0$.
\end{lemma}

\begin{proof}
For finite measures the weak-convergence assertion is the usual Portmanteau theorem, after normalization of the total masses \cite[Chapter 1]{Billingsley1999}. If the limiting mass is zero, the conclusion is immediate. For the vague case, choose
$\chi\in C_c$ with $0\le\chi\le1$ and $\chi=1$ on the support of $G$. Then $\chi m_j\to\chi m$ weakly, since $f\chi\in C_c$ for every bounded continuous $f$. Apply the finite-measure result to $G$ and use $\chi G=G$. Taking $G=\mathbbm{1}_B$ gives the last assertion.
\end{proof}
We use the Fourier convention
\[
\widehat\psi(t)=\int_\R e^{-its}\psi(s)\,ds,
\qquad
\psi(s)=\frac1{2\pi}\int_\R e^{its}\widehat\psi(t)\,dt.
\]
\begin{proof}[Proof of Theorem~\ref{thm:renewal-limit}]
The proof is split into three parts. First, we prove convergence of $W_\mu^x$ to $\nu_0\otimes \Leb$ against tensor-product test functions $\varphi$ and $\psi$. Then we use a result that allows us to transfer this convergence to arbitrary $F\in C_c(\mathbf P(\R^d)\times\R)$. Finally, we remove the weighting $e^{q_\mu s}$ to achieve vague convergence of $|\lambda_\mu| U_\mu^x$ to $\nu_0\otimes \Leb$.

Take $\varphi\in\mathcal C^\alpha(\mathbf P(\R^d))$ and
$\psi\in\mathcal S(\R)$ with $\widehat\psi\in C_c^\infty(\R)$ (where $\mathcal S(\R)$ denotes the set of Schwartz functions on $\R$).
Let $r$ be the radius in Lemma~\ref{lem:critical-resolvent}, and choose
$T>r$ with $\supp\widehat\psi\subset[-T,T]$.
Definitions \eqref{eq:weighted-renewal-measure} and
\eqref{eq:renewal-measure} give
\begin{align*}
\int\varphi([v])\psi(s)\,dW_\mu^x
&=|\lambda_\mu|\int e^{q_\mu s}\varphi([v])\psi(s)\,dU_\mu^x\\
&=|\lambda_\mu|\sum_{n\ge0}\E_\mu
 [e^{q_\mu S_n(x)}\varphi(\theta_n)\psi(S_n(x))].
\end{align*}
Insert Fourier inversion at the random argument,
\[
\psi(S_n(x))
=\frac1{2\pi}\int_{-T}^T
\widehat\psi(t)e^{itS_n(x)}\,dt.
\]
The interchange of sum, expectation and integral is justified by
\eqref{eq:weight-negative-drift}, since
\[
\sum_{n\ge0}\int_{-T}^T |\widehat\psi(t)|\,
\E_\mu[e^{q_\mu S_n(x)}|\varphi(\theta_n)|]\,dt
\le \|\widehat\psi\|_1\|\varphi\|_\infty
\sum_{n\ge0}\E_\mu[e^{q_\mu S_n(x)}]<\infty.
\]
We obtain
\begin{align*}
\int\varphi\psi\,dW_\mu^x
&=\frac{|\lambda_\mu|}{2\pi}
\int_{-T}^T\widehat\psi(t)
\sum_{n\ge0}\E_\mu[e^{(q_\mu+it)S_n(x)}\varphi(\theta_n)]\,dt\\
&=\frac{|\lambda_\mu|}{2\pi}
\int_{-T}^T\widehat\psi(t)
\sum_{n\ge0}(\PP_{q_\mu+it,\mu}^{\,n}\varphi)([x])\,dt,
\end{align*}
where the second equality is the operator identity
\eqref{eq:operator-Pn}. For this fixed $T$, Lemmas~\ref{lem:away-frequency}
and \ref{lem:critical-resolvent}, together with $q_\mu\to0$, give
on a sufficiently small law neighbourhood
\[
\sum_{n\ge0}\PP_{q_\mu+it,\mu}^{\,n}
=(I-\PP_{q_\mu+it,\mu})^{-1},\qquad |t|\le T.
\]
Thus
\begin{equation}\label{eq:weighted-Fourier}
\int\varphi\psi\,dW_\mu^x
=\frac{|\lambda_\mu|}{2\pi}
\int_{-T}^T\widehat\psi(t)
\bigl((I-\PP_{q_\mu+it,\mu})^{-1}\varphi\bigr)([x])\,dt.
\end{equation}

Split the integral in \eqref{eq:weighted-Fourier} into
$|t|\le r$ and $r<|t|\le T$, denoting the two contributions by $I_{\mathrm{near},\mu}^x$ and
$I_{\mathrm{away},\mu}^x$.
By \eqref{eq:away-resolvent} in Lemma~\ref{lem:away-frequency},
\[
|I_{\mathrm{away},\mu}^x|
\le C_T|\lambda_\mu|\,\|\varphi\|_{\mathcal C^\alpha}\|\widehat\psi\|_1
\longrightarrow0.
\]
The nonzero limit therefore comes from frequencies near zero.
In the near integral substitute $t=|\lambda_\mu|u$ and put
\[
B_\mu(u):=|\lambda_\mu|^2
(I-\PP_{q_\mu+i|\lambda_\mu|u,\mu})^{-1}.
\]
Then
\begin{equation}\label{eq:weighted-small-frequency}
I_{\mathrm{near},\mu}^x
=\frac1{2\pi}\int_{|u|\le r/|\lambda_\mu|}
\widehat\psi(|\lambda_\mu|u)
(B_\mu(u)\varphi)([x])\,du.
\end{equation}
The limit \eqref{eq:critical-resolvent-limit} says exactly that
\begin{equation}\label{eq:renewal-rescaled-operator}
B_\mu(u)\longrightarrow
\frac{2\sigma_0^2}{1+\sigma_0^4u^2}\Pi_{0,\mu_0}
\quad\text{in operator norm, uniformly for bounded }u.
\end{equation}
Moreover, substituting $t=|\lambda_\mu|u$ into
\eqref{eq:resolvent-denominator} gives
\begin{equation}\label{eq:renewal-u-majorant}
\|B_\mu(u)\|_{\mathcal C^\alpha}\le\frac{C}{1+u^2},
\qquad |u|\le r/|\lambda_\mu|.
\end{equation}

Define
\[
f_\mu^x(u):=\frac1{2\pi}
\mathbbm{1}_{\{|u|\le r/|\lambda_\mu|\}}
\widehat\psi(|\lambda_\mu|u)(B_\mu(u)\varphi)([x]),
\]
with value zero outside the integration range in
\eqref{eq:weighted-small-frequency}. Thus
$I_{\mathrm{near},\mu}^x=\int_\R f_\mu^x(u)\,du$. Since
$\widehat\psi$ is continuous and
$\Pi_{0,\mu_0}\varphi=\nu_0(\varphi)\mathbbm{1}_{\mathbf{P}(\mathbb{R}^d)}$ by
Proposition~\ref{prop:analytic-spectral-data},
\eqref{eq:renewal-rescaled-operator} yields
\[
f_\mu^x(u)\longrightarrow
f_0(u):=\frac{\widehat\psi(0)\nu_0(\varphi)}{2\pi}
\frac{2\sigma_0^2}{1+\sigma_0^4u^2}
\]
uniformly on each bounded $u$-interval. Both integrands are bounded
by $C_{\varphi,\psi}/(1+u^2)$, by
\eqref{eq:renewal-u-majorant} and
$\|\widehat\psi\|_\infty\le\|\psi\|_1$. For fixed $M>0$,
\[
\sup_x\left|I_{\mathrm{near},\mu}^x-\int_\R f_0(u)\,du\right|
\le 2M\sup_{\substack{|u|\le M\\x}}|f_\mu^x(u)-f_0(u)|
+C\int_{|u|>M}\frac{du}{1+u^2}.
\]

Using a change of variable $v=\sigma_0^2 u$, \begin{equation}\label{eq:kernel-mass-one}
\frac1{2\pi}\int_\R\frac{2\sigma_0^2}{1+\sigma_0^4u^2}\,du
=\frac1\pi\int_\R\frac{dv}{1+v^2}=1.
\end{equation}

Hence, first taking $\mu\to\mu_0$ with $M$ fixed and then $M\to\infty$, the away integral vanishes, and since $\widehat\psi(0)=\int\psi(s)\,ds$, we end up with
\begin{equation}\label{eq:weighted-tensor-limit}
\int\varphi([v])\psi(s)\,dW_\mu^x
\longrightarrow\nu_0(\varphi)\int_\R\psi(s)\,ds.
\end{equation}
Thus the limiting projection supplies $\nu_0$, and the Fourier
transform at zero supplies Lebesgue measure.

Set $X=\mathbf P(\R^d)$. Consider
\[
\mathcal H=\{\psi\in\mathcal S(\R):\widehat\psi\in C_c^\infty(\R)\},
\qquad
\mathcal A=\operatorname{span}\{\varphi(v)\psi(s):
\varphi\in\mathcal C^\alpha(X),\ \psi\in\mathcal H\}.
\]
We have proved convergence against this test class. We now need to apply a transfer theorem by L\"ohr and Rippl \cite[Theorem 2.3 and Remark 2.2]{LoehrRippl2016}.
Set $E=X\times\R$, equipped with
$d_E((v,s),(w,t))=d_{\mathbf P}(v,w)+|s-t|$. Since $X$ is compact,
$E$ is separable and its closed bounded sets are compact (this establishes the equivalence between vague convergence and their weak-$\#$ convergence by the remark).
We verify the three hypotheses of this result. The algebra
$\mathcal A\subset C_0(E)$ is closed under multiplication and complex
conjugation, separates points and vanishes nowhere. By the
Stone--Weierstrass theorem it is uniformly dense in $C_0(E)$ and hence
induces the topology of $E$. Every bounded subset of $E$ lies in a strip
$X\times[-R,R]$. Choose $\psi\in\mathcal H$ with $\psi(0)\ne0$;
a suitable dilation is bounded away from zero on $[-R,R]$.
The corresponding function $\mathbbm{1}_{\mathbf{P}(\mathbb{R}^d)}\otimes\psi\in\mathcal A$ therefore
satisfies the third hypothesis.

Each $W_\mu^x$ is a positive finite measure by
\eqref{eq:weighted-total-mass}, and $\nu_0\otimes\Leb$ is finite
on bounded sets. Every element of $\mathcal A$ is integrable against
$W_\mu^x$, since this measure is finite, and against $\nu_0\otimes\Leb$,
since its real-coordinate factors are Schwartz functions.
Equation \eqref{eq:weighted-tensor-limit} gives convergence on $\mathcal A$. The cited theorem therefore yields
vague convergence. Applying it along arbitrary sequences
$\mu_j\to\mu_0$ with $\lambda_{\mu_j}\ne0$ and
$x_j\in\mathbb S^{d-1}$ gives
\begin{equation}\label{eq:localized-Fourier-limit}
\lim_{\substack{\mu\to\mu_0\\\lambda_\mu\ne0}}\sup_x
\left|\int F\,dW_\mu^x-\int F\,d(\nu_0\otimes\Leb)\right|=0,
\qquad F\in C_c(E).
\end{equation}

Final step: removing the weight. Fix $R>0$ and choose a nonnegative $\chi_R\in C_c(\R)$ with
$\chi_R\ge1$ on $[-R,R]$. Applying \eqref{eq:localized-Fourier-limit}
to $F(v,s)=\chi_R(s)$ gives, for all sufficiently close noncritical $\mu$,
\begin{equation}\label{eq:weighted-local-mass}
\sup_x W_\mu^x(X\times[-R,R])
\le\sup_x\int\chi_R(s)\,dW_\mu^x
\le 1+\int_\R\chi_R(s)\,ds=:C_R.
\end{equation}
The identity $|\lambda_\mu|\,dU_\mu^x=e^{-q_\mu s}\,dW_\mu^x$
from \eqref{eq:weighted-renewal-measure} also gives
\begin{equation}\label{eq:uniform-local-mass}
\sup_x |\lambda_\mu|U_\mu^x(X\times[-R,R])
\le e^{|q_\mu|R}C_R\le C'_R<\infty,
\end{equation} since $q_\mu\to 0$.
For $F\in C_c(E)$ supported in $|s|\le R$, the same identity and
\eqref{eq:weighted-local-mass} yield
\begin{align*}
\left||\lambda_\mu|\int F\,dU_\mu^x
-\int F\,dW_\mu^x\right|
&=\left|\int(e^{-q_\mu s}-1)F\,dW_\mu^x\right|\\
&\le C_R\|F\|_\infty
\sup_{|s|\le R}|e^{-q_\mu s}-1|\longrightarrow0.
\end{align*}
Together with \eqref{eq:localized-Fourier-limit}, this proves
\eqref{eq:renewal-limit}.
\end{proof}

\section{Escape Counts}\label{sec:escape-mean}

In the last section we proved the vague convergence $|\lambda_\mu|U_\mu^{x} \to \nu_0\otimes\Leb$. We now need to define a set that encodes an escape at time $n$. However, since an escape would depend on the choice of matrix at time $n+1$, we need to extend the renewal measure to encode the choice of the next matrix. In this context, we define the set $\cC$ below that entails a nonpositive logarithmic norm at time $n$ and a positive log norm at time $n+1$. We will prove in this section that the extended measure of this set will correspond exactly to $|\lambda_\mu| \mathbb{E}_\mu [N_x]$, and the vague convergence theorem will yield convergence to a constant $\kappa_0$, hence proving the exact critical exponent $1$ for this observable as the system approaches the critical measure $\mu_0$. We will also establish the exponential scaling limit of this observable.

Let
\[
\mathsf E:=\mathbf P(\R^d)\times\R\times S,
\qquad M_0:=\nu_0\otimes\Leb\otimes\mu_0.
\]
For every bounded Borel function $G$ that is $M_0$-integrable, write
$M_0(G):=\int_{\mathsf E}G\,dM_0$.
For a nearby noncritical law $\mu$ and $x\in\mathbb S^{d-1}$, define
\begin{equation}\label{eq:extended-measure}
\widetilde U_\mu^x:=|\lambda_\mu|U_\mu^x\otimes\mu.
\end{equation}
This is a Radon measure on $\mathsf E$. The additional coordinate
records the next matrix, which is independent of the current state.

\begin{lemma}[Vague Convergence of the extended measures]\label{lem:extended-measures}
For every $G\in C_c(\mathsf E)$ and compact $J\subset\R$, write
$G_a([v],s,A):=G([v],s+a,A)$. Then
\begin{equation}\label{eq:extended-vague}
\lim_{\substack{\mu\to\mu_0\\\lambda_\mu\ne0}}
\sup_{\substack{x\in\mathbb S^{d-1}\\a\in J}}
\left|\int G_a\,d\widetilde U_\mu^x-\int G\,dM_0\right|=0.
\end{equation}
In particular, $\widetilde U_\mu^x\to M_0$ vaguely, uniformly in $x$.
\end{lemma}

\begin{proof}
Define
\[
F_\mu([v],s):=\int_S G([v],s,A)\,d\mu(A),
\qquad
F_0([v],s):=\int_S G([v],s,A)\,d\mu_0(A).
\]
These functions are continuous and vanish outside one fixed compact
strip since $G\in C_c(\mathsf E)$. Fubini's theorem gives
\[
\int G_a\,d\widetilde U_\mu^x
=|\lambda_\mu|\int F_\mu([v],s+a)\,dU_\mu^x,
\qquad
\int G\,dM_0=\int F_0\,d(\nu_0\otimes\Leb).
\]
Split the difference as
\begin{align*}
\int G_a\,d\widetilde U_\mu^x-\int G\,dM_0
&=|\lambda_\mu|\int(F_\mu-F_0)([v],s+a)\,dU_\mu^x\\
&\quad+\left(|\lambda_\mu|\int F_0([v],s+a)\,dU_\mu^x
-\int F_0\,d(\nu_0\otimes\Leb)\right).
\end{align*}
The second term tends to zero uniformly in $x$ and $a\in J$ by Theorem \ref{thm:renewal-limit}.
For the first term, note that \[
\|F_\mu-F_0\|_\infty
\le\|G\|_\infty\|\mu-\mu_0\|_{\TV}\longrightarrow0.
\]
The functions $F_\mu$ and $F_0$ are compactly supported. The local bound \eqref{eq:uniform-local-mass} therefore
gives, for all sufficiently nearby noncritical laws, \[
\sup_{x\in \mathbb{S}^{d-1},a\in J}
\left||\lambda_\mu|\int(F_\mu-F_0)([v],s+a)\,dU_\mu^x\right|
\le C_{G,J}\|F_\mu-F_0\|_\infty\longrightarrow0.
\]
Together these estimates prove \eqref{eq:extended-vague}.
\end{proof}

We now have the tools to prove the critical exponent $1$ for the expected number of escapes.

\begin{proof}[Proof of Theorem \ref{MainTheoremEscapesCritical} (1)]
    Define the outward-crossing set
\begin{equation}\label{eq:crossing-set}
\cC:=\{([v],s,A)\in\mathsf E:s\le0<s+\xi(A,[v])\}.
\end{equation}
At time $n$, the current state is $(\theta_n,S_n(x))$ and the next
logarithmic norm is $S_n(x)+\xi(A_{n+1},\theta_n)$.
Thus $\cC$ records a crossing of the unit sphere from inside to outside.
By compactness of the set $S$, we have \[ |\xi(A,[v])|\leq \log \eta \quad \text{for all } A\in S, [v] \in \mathbf{P}(\R^d) .\]
This implies that $\cC\subset \mathbf{P}(\R^d) \times [-\log \eta, 0] \times S $, which is a compact set. Since $\xi$ is a continuous function, the boundary of $\cC$ has zero $M_0$ measure since
\begin{equation}\label{eq:crossing-boundary}
\partial\cC\subset\{s=0\}\cup\{s+\xi(A,[v])=0\}.
\end{equation}
The two assumptions of Lemma \ref{lem:vague-portmanteau} are satisfied, so applying this lemma to the vague convergence of the extended measure (Lemma \ref{lem:extended-measures} with $a=0$) gives
\begin{equation}\label{eq:crossing-portmanteau}
\sup_{x\in\mathbb S^{d-1}}
\left|\widetilde U_\mu^x(\cC)-M_0(\cC)\right|
\longrightarrow0.
\end{equation}
We now need to compute each of $\widetilde U_\mu^x(\cC)$ and $M_0(\cC)$. Indeed, the former coincides with the normalized expected number of escapes starting from some $x\in \mathbb{S}^{d-1}$ under the measure $\Prob_\mu$.
Independence of $A_{n+1}$ from the first $n$ matrices gives
\begin{equation}\label{eq:joint-law-current-next}
\operatorname{Law}_\mu(\theta_n,S_n(x),A_{n+1})
=\operatorname{Law}_\mu(\theta_n,S_n(x))\otimes\mu.
\end{equation}
Summing over $n$ and using \eqref{eq:renewal-measure} yields the
identity of locally finite measures
\begin{equation}\label{eq:extended-law-sum}
U_\mu^x\otimes\mu
=\sum_{n=0}^\infty
\operatorname{Law}_\mu(\theta_n,S_n(x),A_{n+1}).
\end{equation}
Membership in $\cC$ is exactly the event $S_n(x)\le0<S_{n+1}(x)$.
Evaluating \eqref{eq:extended-law-sum} on $\cC$ and using
the definition of $N_x$ therefore gives
\begin{equation}\label{eq:crossing-identity}
\widetilde U_\mu^x(\cC)
=|\lambda_\mu|\sum_{n=0}^\infty
\Prob_\mu(S_n(x)\le0<S_{n+1}(x))
=|\lambda_\mu|\E_\mu[N_x].
\end{equation}

For fixed $u\in\R$,
\[
\{s\in\R:s\le0<s+u\}
=\begin{cases}(-u,0],&u>0,\\ \varnothing,&u\le0\end{cases}.
\]
Its Lebesgue measure is $u^+$. Hence
\begin{align}
M_0(\cC)
&=\int_{\mathbf P(\R^d)}\int_S\int_\R
\mathbbm{1}_{\{s\le0<s+\xi(A,[v])\}}\,ds\,d\mu_0(A)\,d\nu_0([v])\notag\\
&=\int_{\mathbf P(\R^d)}\int_S\xi(A,[v])^+
\,d\mu_0(A)\,d\nu_0([v])=:\kappa_0.
\label{eq:crossing-flux}
\end{align}
Furstenberg's formula for the top Lyapunov exponent and the fact that $\lambda_{\mu_0}=0$ give
\[
\int_{\mathbf P(\R^d)}\int_S\xi(A,[v])\,d\mu_0(A)\,d\nu_0([v])= \lambda_{\mu_0}=0.
\]
The positive and negative parts of $\xi$ therefore have equal integrals,
so
\[
\kappa_0=\frac12\int_{\mathbf P(\R^d)}\int_S
|\xi(A,[v])|\,d\mu_0(A)\,d\nu_0([v]),
\]
Finally, using \eqref{eq:crossing-portmanteau}, we obtain
\[ \sup_{x\in \mathbb{S}^{d-1}} \left||\lambda_\mu| \E_{\mu}[N_x] - \kappa_0\right| \to 0 \quad \text{as } \mu \to \mu_0  \]
along noncritical measures.
Since $|\xi(A,[v])|\leq \log \eta$, we have $\kappa_0<\infty$. We also have $\kappa_0>0$: otherwise, $\xi=0$ for
$\mu_0\otimes\nu_0$-almost every pair. With a stationary initial direction,
every cocycle increment would then be zero almost surely, contradicting
the nondegenerate central limit theorem.

It remains to show that the number of escapes is almost surely infinite at $\mu_0$. This is a direct consequence of the law of iterated logarithms for random matrix products.
By \cite[Chapter~14]{benoist2016random}, for all $x\in \mathbb{S}^{d-1}$, the set of accumulation points of the sequence \[ \frac{\log \|A_\omega^n x\| - \lambda_{\mu_0} n}{\sqrt{2n \log \log n}}  \]
is the interval $[-\sigma_{\mu_0},\sigma_{\mu_0}]$ $\Prob_{\mu_0}$-a.s. Since $\sigma_{\mu_0}>0$ and $\lambda_{\mu_0} =0$, we have \[ \Prob_{\mu_0} \left[ \liminf_{n\to \infty}  \|A_\omega^n x\|<1 \text{ and } \limsup_{n\to \infty}  \|A_\omega^n x\|>1  \right]=1. \]
Hence, $\Prob_{\mu_0} [N_x=\infty]=1$.
\end{proof}

We now move to the exponential scaling limit.

\begin{proof}[Proof of Theorem \ref{MainTheoremEscapesCritical} (2)]

For $y\ne0$, define $N_y$ by the same crossing formula with initial
vector $y$. Let $\mathcal{A}_\eta := \{ y\in \R^d: 1\le \|y\| \le \eta \}$. Every point reached immediately after an escape lies in this strip since $\|A\|\leq \eta$ for all $A\in S$. For $y\in \mathcal{A}_\eta$, we have $\log \|A_\omega^n y\| = \log \|y\| + S_n(y/\|y\|)$. Using the same set $\cC$ defined above, we have
\[ |\lambda_\mu| \mathbb{E}_\mu [N_y] = \int_{\mathsf{E}} \mathbbm{1}_{\cC} ([v], s+\log\|y\|,A) d\widetilde U_\mu^{y/\|y\|} ([v],s,A). \]
Lemma \ref{lem:extended-measures} now implies
\begin{equation}\label{eq:shifted-escape-mean}
\sup_{y\in\mathcal A_\eta}
\bigl||\lambda_\mu| \E_\mu [N_y] -\kappa_0\bigr|
\longrightarrow0.
\end{equation}
Set $H_\mu(y):=\mathbb E_\mu[N_y]$, where the expectation is taken
for a new product started at the fixed vector $y$. Let
$m_\mu^+ = \sup_{y\in \mathcal{A}_\eta} H_\mu(y)$ and
$m_\mu^- = \inf_{y\in \mathcal{A}_\eta} H_\mu(y)$. By \eqref{eq:shifted-escape-mean}, $|\lambda_\mu|m_\mu^{\pm} \to \kappa_0>0$, so for all nearby noncritical measures $\mu$, we have $0<m_\mu^- \le m_\mu^+ <\infty$.

For $k\ge 1$, let $T_k$ be the time right after the $k^{th}$ escape (a stopping time). We have \[ \{T_k<\infty\} = \{N_x\ge k\}. \]
Note that on $\{T_k<\infty\}$ we have $1\leq \|A_\omega^{T_k} x\|\le \eta$. Let $\mathcal{F}_n = \sigma(A_1,\dots, A_n)$. Define \[ R_k = \E_\mu [(N_x-k)^+] \]
as the expected number of escapes happening after the $k^{th}$ escape, and note that $R_0 = \E_\mu [N_x]$.

On $\{T_k<\infty\}$, exactly $k$ escapes have occurred by time $T_k$. Conditionally on $\mathcal F_{T_k}$, the remaining matrices are independent with law
$\mu$, so the remaining $N_x-k$ escapes have the law of the escape
count for a new product started at $A_\omega^{T_k} x$. The strong Markov property and the tower property yield
\begin{align}
R_k
&=\E_\mu\!\left[\mathbbm{1}_{\{T_k<\infty\}}(N_x-k)\right]\notag\\
&=\E_\mu\!\left[
\E_\mu\!\left[\mathbbm{1}_{\{T_k<\infty\}}(N_x-k)
\mid\mathcal F_{T_k}\right]\right]\notag\\
&=\E_\mu\!\left[
\mathbbm{1}_{\{T_k<\infty\}}\E_\mu[N_x-k\mid\mathcal F_{T_k}]
\right]\notag\\
&=\E_\mu\!\left[\mathbbm{1}_{\{T_k<\infty\}} H_\mu(A_\omega^{T_k}x) \right].
\label{eq:escape-residual-mean}
\end{align}
Here the integrand is set to zero on $\{T_k=\infty\}$. This implies
\begin{equation}\label{eq:escape-residual-bounds}
    m_\mu^- \Prob_\mu [N_x\ge k] \le R_k \le  m_\mu^+ \Prob_\mu [N_x\ge k],
\end{equation}
for all $k\geq 1$. Since $R_{k-1} - R_k = \E_\mu[\mathbbm{1}_{T_k<\infty}] = \Prob_\mu [N_x\geq k]$, we obtain the recursive bounds
\[
r_\mu^-R_{k-1}\le R_k\le r_\mu^+R_{k-1},
\qquad \text{where }r_\mu^\pm:=\frac{m_\mu^\pm}{1+m_\mu^\pm}.
\]
Equation \eqref{eq:escape-residual-bounds} also yields $R_k/m_\mu^+ \leq \Prob_\mu [N_x\geq k] \le R_k/m_\mu^-$. Using these bounds and iterating the above starting $R_0 = \E_\mu[N_x]$ gives
\begin{equation}\label{eq:escape-geometric-bounds}
\frac{\E_\mu[N_x]}{m_\mu^+}(r_\mu^-)^k
\le\Prob_\mu(N_x\ge k)
\le\frac{\E_\mu[N_x]}{m_\mu^-}(r_\mu^+)^k,
\qquad k\ge1.
\end{equation}
For $t\geq 0$, set $k_\mu(t) := \lfloor t/|\lambda_\mu| \rfloor +1 $ so that \[ \Prob_\mu[|\lambda_\mu| N_x > t] = \Prob_\mu [N_x\geq k_\mu(t)] \]
because $N_x$ is integer-valued. Moreover, $|\lambda_\mu| k_\mu(t) \to t$ as $\mu \to \mu_0$. When taking this limit, \eqref{eq:shifted-escape-mean} implies
$\E_\mu[N_x]/m_\mu^\pm\to1$ uniformly in $x\in \mathbb{S}^{d-1}$. Moreover,
\[
\frac{\log r_\mu^\pm}{|\lambda_\mu|}
=-\frac{\log(1+1/m_\mu^\pm)}{|\lambda_\mu|} = -\frac{1}{|\lambda_\mu|m_\mu^\pm } \frac{\log(1+1/m_\mu^\pm)}{1/m_\mu^\pm}
\longrightarrow-\frac1{\kappa_0}
\]
using $\log(1+z)\sim z$ as $z\to 0$. Both bounds in \eqref{eq:escape-geometric-bounds} therefore converge to $e^{-t/\kappa_0}$. We conclude that
\begin{equation}\label{eq:escape-tail-limit}
\sup_{x\in\mathbb S^{d-1}}
\left|\Prob_\mu(|\lambda_\mu|N_x>t)-e^{-t/\kappa_0}\right|
\longrightarrow0,
\qquad t\ge0.
\end{equation}
This proves the assertion \[|\lambda_\mu| N_x \xrightarrow{\mathrm{d}} \kappa_0 \mathcal{E}\] where $\mathcal{E}\sim\operatorname{Exp}(1)$, as $\mu \to \mu_0$ along noncritical measures.

\end{proof}

\section{Time Spent in a Fixed Strip}\label{sec:strip-visits}

We now prove Theorem \ref{MainTheoremTimeSpentInStripCritical}. As in
Section \ref{sec:escape-mean}, the critical exponent follows by
evaluating the renewal measure on a suitable set. Here the event of
interest depends only on the current logarithmic norm, so we can apply
Theorem \ref{thm:renewal-limit} directly. The exponential limit then
follows from the same restart argument used for the number of escapes.

Fix $a,b\in\R$ with $a<b$, and write
\[
\mathcal A_{a,b}:=\{y\in\R^d:e^a\leq\|y\|\leq e^b\}.
\]
For $y\ne0$, define
\[
L_y(a,b):=\sum_{n=0}^{\infty}
\mathbbm{1}_{\{e^a\leq\|A_\omega^n y\|\leq e^b\}}.
\]
This counts the visits to $\mathcal A_{a,b}$ by the product started
at $y$, with time zero included when $y\in\mathcal A_{a,b}$.

\begin{proof}[Proof of Theorem \ref{MainTheoremTimeSpentInStripCritical} (1)]
Define
\[
\mathcal D_{a,b}:=\mathbf P(\R^d)\times[a,b].
\]
By \eqref{eq:renewal-measure} and Tonelli's theorem,
\begin{equation}\label{eq:strip-occupation-identity}
U_\mu^x(\mathcal D_{a,b})
=\sum_{n=0}^{\infty}\Prob_\mu(a\leq S_n(x)\leq b)
=\E_\mu[L_x(a,b)].
\end{equation}
The set $\mathcal D_{a,b}$ is compact and its boundary is
$\mathbf P(\R^d)\times\{a,b\}$, which has zero
$\nu_0\otimes\Leb$ measure. Applying
Lemma \ref{lem:vague-portmanteau} to the convergence in
Theorem \ref{thm:renewal-limit} therefore gives
\begin{equation}\label{eq:strip-mean}
\sup_{x\in\mathbb S^{d-1}}
\left||\lambda_\mu|\E_\mu[L_x(a,b)]-(b-a)\right|
\longrightarrow0
\end{equation}
as $\mu\to\mu_0$ along noncritical measures, since
$(\nu_0\otimes\Leb)(\mathcal D_{a,b})=b-a$.
This proves \eqref{Expected_Strip_Time_Critical_Exponent}.

For the exponential limit, we will also need the same estimate when the starting point lies anywhere in $\mathcal A_{a,b}$. Write $y=e^c v$, where $v\in\mathbb S^{d-1}$ and $c\in[a,b]$.
Since $\log\|A_\omega^n y\|=c+S_n(v)$,
\[
\E_\mu[L_y(a,b)]
=U_\mu^v\bigl(\mathbf P(\R^d)\times[a-c,b-c]\bigr).
\]
For each fixed $c$, the preceding argument applies to this translated
interval, whose length is still $b-a$. The convergence is uniform
for $c\in[a,b]$ (just like for escapes):
\begin{equation}\label{eq:shifted-strip-mean}
\sup_{y\in\mathcal A_{a,b}}
\left||\lambda_\mu|\E_\mu[L_y(a,b)]-(b-a)\right|
\longrightarrow0.
\end{equation}

It remains to prove that the number of visits is almost surely
infinite at $\mu_0$. The local limit theorem
\cite[Theorem 2.4]{XGL_Mod}, applied at $\lambda_{\mu_0}=0$, gives
\[
\sup_{\substack{v\in\mathbb S^{d-1}\\0\leq c\leq\log\eta}}
\left|
\sqrt n\,\Prob_{\mu_0}(a\leq c+S_n(v)\leq b)
-\frac{b-a}{\sigma_0\sqrt{2\pi}}
\right|\longrightarrow0.
\]
In particular, there exist an integer $m\geq1$ and $p>0$ such that
\begin{equation}\label{eq:critical-strip-hitting}
\inf_{1\leq\|y\|\leq\eta}
\Prob_{\mu_0}(A_\omega^m y\in\mathcal A_{a,b})\geq p.
\end{equation}
By Theorem \ref{MainTheoremEscapesCritical} (1), there are almost surely infinitely many escapes. Every point immediately after an escape has norm in $(1,\eta]$, so
$1\leq\|A_\omega^n x\|\leq\eta$ for infinitely many $n$, almost surely.
Let
\[
E_n:=\{A_\omega^n x\in\mathcal A_{a,b}\},
\qquad \mathcal F_n:=\sigma(A_1,\dots,A_n).
\]
Conditionally on $\mathcal F_n$, the next $m$ matrices are independent
with law $\mu_0$. Thus \eqref{eq:critical-strip-hitting} gives
\[
\Prob_{\mu_0}(E_{n+m}\mid\mathcal F_n)
\geq p\,\mathbbm{1}_{\{1\leq\|A_\omega^n x\|\leq\eta\}},
\]
and consequently
\begin{equation}\label{eq:critical-strip-conditional-sum}
\sum_{n=0}^{\infty}
\Prob_{\mu_0}(E_{n+m}\mid\mathcal F_n)=\infty
\qquad\Prob_{\mu_0}\text{-a.s.}
\end{equation}
Apply the conditional Borel--Cantelli lemma
\cite[Theorem~4.5.5]{Durrett2019} along each residue class
modulo $m$: for $r=0,\dots,m-1$, use the events
$E_{r+jm}$, $j\geq1$, and the filtration
$(\mathcal F_{r+jm})_{j\geq0}$. On the event
\[
\sum_{j=0}^{\infty}
\Prob_{\mu_0}(E_{r+(j+1)m}\mid\mathcal F_{r+jm})=\infty,
\]
infinitely many of these events occur, almost surely. By
\eqref{eq:critical-strip-conditional-sum}, at least one of the
finitely many residue classes has a divergent sum, almost surely.
It follows that $E_n$ occurs infinitely often, proving
$\Prob_{\mu_0}[L_x(a,b)=\infty]=1$.
\end{proof}

\begin{proof}[Proof of Theorem \ref{MainTheoremTimeSpentInStripCritical} (2)]
We use the same strategy developed for the exponential limit of escape counts. For $y\in\mathcal A_{a,b}$, set
\[
H_\mu(y):=\E_\mu[L_y(a,b)]-1,
\qquad
m_\mu^-:=\inf_{y\in\mathcal A_{a,b}}H_\mu(y),
\qquad
m_\mu^+:=\sup_{y\in\mathcal A_{a,b}}H_\mu(y).
\]
The subtraction of $1$ removes the visit at time zero. By
\eqref{eq:shifted-strip-mean},
$|\lambda_\mu|m_\mu^\pm\to b-a>0$, so
$0<m_\mu^-\leq m_\mu^+<\infty$ for all sufficiently nearby
noncritical measures.

For $k\geq1$, let $\tau_k$ be the time of the $k$th visit to
$\mathcal A_{a,b}$, and write
\[
R_k:=\E_\mu[(L_x(a,b)-k)^+],\qquad k\geq0.
\]
On $\{\tau_k<\infty\}=\{L_x(a,b)\geq k\}$, the point
$A_\omega^{\tau_k}x$ lies in $\mathcal A_{a,b}$. Restarting there
counts the current visit again, whereas $L_x(a,b)-k$ counts only
subsequent visits. Thus the strong Markov property gives
\begin{equation}\label{eq:strip-residual-mean}
R_k
=\E_\mu\!\left[
\mathbbm{1}_{\{\tau_k<\infty\}}
H_\mu(A_\omega^{\tau_k}x)
\right],\qquad k\geq1.
\end{equation}
Here the integrand is set to zero on $\{\tau_k=\infty\}$.
Consequently,
\[
\begin{aligned}
m_\mu^-\Prob_\mu(L_x(a,b)\geq k)
&\leq R_k
\leq m_\mu^+\Prob_\mu(L_x(a,b)\geq k),\\
R_{k-1}-R_k&=\Prob_\mu(L_x(a,b)\geq k).
\end{aligned}
\]
These are the same relations used to obtain \eqref{eq:escape-geometric-bounds} from \eqref{eq:escape-residual-bounds}. Applying that calculation with
$R_0=\E_\mu[L_x(a,b)]$ and $r_\mu^\pm=m_\mu^\pm/(1+m_\mu^\pm)$ yields
\begin{equation}\label{eq:strip-geometric-bounds}
\frac{\E_\mu[L_x(a,b)]}{m_\mu^+}(r_\mu^-)^k
\leq\Prob_\mu(L_x(a,b)\geq k)
\leq\frac{\E_\mu[L_x(a,b)]}{m_\mu^-}(r_\mu^+)^k,
\qquad k\geq1.
\end{equation}
For $t\geq0$, take
$k_\mu(t)=\lfloor t/|\lambda_\mu|\rfloor+1$.
Equations \eqref{eq:strip-mean} and \eqref{eq:shifted-strip-mean}
give $\E_\mu[L_x(a,b)]/m_\mu^\pm\to1$ uniformly in $x$.
The same logarithmic calculation used for the escape count gives
\[
k_\mu(t)\log r_\mu^\pm\longrightarrow-\frac{t}{b-a}.
\]
Since
$\{|\lambda_\mu|L_x(a,b)>t\}=\{L_x(a,b)\geq k_\mu(t)\}$, both bounds in \eqref{eq:strip-geometric-bounds} converge to $e^{-t/(b-a)}$. Hence
\begin{equation}\label{eq:strip-tail-limit}
\sup_{x\in\mathbb S^{d-1}}
\left|\Prob_\mu(|\lambda_\mu|L_x(a,b)>t)-e^{-t/(b-a)}\right|
\longrightarrow0,\qquad t\geq0.
\end{equation}
This proves $|\lambda_\mu|L_x(a,b)\xrightarrow{\mathrm d}
(b-a)\mathcal E$, where $\mathcal E$ is exponential with mean $1$,
uniformly in $x\in\mathbb S^{d-1}$.
\end{proof}

\section{Opposite-phase residence}\label{sec:residence}

This section is dedicated to the proof of Theorem \ref{opposite-exponents} which establishes the critical exponent and scaling limit of the opposite-phase residence 
\begin{equation}\label{eq:opposite-observables}
\Topp:=\sum_{n=0}^{\infty}
\mathbbm{1}_{\{-\operatorname{sgn}(\lambda_\mu) S_n(x)>0\}},
\end{equation}
for $x\in \mathbb{S}^{d-1}$ and $\lambda_\mu \ne 0$.
For the process $Z_\mu^x$ defined in \eqref{eq:diffusive-cocycle-process}, notice that \[|\lambda_\mu|^2 \Topp = \int_0^\infty \mathbbm{1}_{ \{ Z_{\mu}^x(t)>0 \} } dt \]

since $-\lambda_\mu S_n(x)=|\lambda_\mu|[-\operatorname{sgn}(\lambda_\mu)S_n(x)]$ and each matrix step occupies a rescaled time interval of length $|\lambda_\mu|^2$. We first prove Theorem~\ref{lem:diffusive-cocycle-limit}, the path limit for $Z_\mu^x$ needed to pass to occupation times. The argument uses the uniform spectral data of Proposition~\ref{prop:analytic-spectral-data} which allows the law $\mu$ and the starting direction to vary together.

\begin{proof}[Proof of Theorem~\ref{lem:diffusive-cocycle-limit}]
    By \cite[Theorem 8.3]{Sagitov2020}, we need to verify three conditions: convergence of the finite dimensional distributions, continuity of the limiting process at $K$, and a tightness criterion of the form
    \begin{equation}\label{tightness_Cond3}
\E_\mu\!\left[|Z_\mu^x(s)-Z_\mu^x(r)|^2|Z_\mu^x(t)-Z_\mu^x(s)|^2\right]
\le 4C_K^2(t-r)^2,
\end{equation}
for any $0\le r \le s \le t \le K$.

We begin with convergence of the finite dimensional distributions. Proposition~\ref{prop:analytic-spectral-data} gives, uniformly for nearby laws,
\begin{equation}
\begin{gathered}\label{diffusive-spectral-input}
\PP_{z,\mu}^{\,m}
=e^{m\Lambda_\mu(z)}\Pi_{z,\mu}+\mathcal Q_{z,\mu}^{\,m},
\qquad \|\mathcal Q_{z,\mu}^{\,m} \|_{\mathcal C^\alpha}\le C\rho^m,
\quad 0<\rho<1,\\
\Lambda_\mu(z)=\lambda_\mu z+\tfrac12\sigma_\mu^2z^2+O(|z|^3),
\qquad \sigma_\mu^2\longrightarrow\sigma_0^2.
\end{gathered}
\end{equation}
For fixed $u\in\R$ and $m=m_\mu$ with
$m|\lambda_\mu|^2\to t>0$, substitution of $z=-iu\lambda_\mu$ gives
\[
m\Lambda_\mu(-iu\lambda_\mu)
=m|\lambda_\mu|^2
\bigl(-iu-\tfrac12\sigma_\mu^2u^2+O(|\lambda_\mu|)\bigr)
\longrightarrow t\bigl(-iu-\tfrac12\sigma_0^2u^2\bigr)
\]
as $\mu \to \mu_0$. Since $\Pi_{-iu\lambda_\mu,\mu}\to\Pi_{0,\mu_0}$ in operator norm,
\begin{equation}\label{eq:diffusive-block-limit}
\PP_{-iu\lambda_\mu,\mu}^{\,m}
\longrightarrow e^{t(-iu-\sigma_0^2u^2/2)}\Pi_{0,\mu_0}
\quad\text{in operator norm},
\end{equation}
noting that $\Pi_{0,\mu_0} \mathbbm{1}_{\mathbf{P}(\mathbb{R}^d)} = \mathbbm{1}_{\mathbf{P}(\mathbb{R}^d)}$.

Now fix $0=t_0<t_1<\dots <t_\ell\le K$. Setting $Y_{t_i} = \sigma_0 B_{t_i}-t_i$, we want to show that \[  (Z_\mu^x(t_1),Z_\mu^x(t_2),\dots,Z_\mu^x(t_\ell)) \xrightarrow{d} (Y_{t_1}, Y_{t_2}, \dots, Y_{t_\ell}  ) .  \]
We achieve this by proving convergence of the joint characteristic functions of the increments. Set $n_j=\lfloor t_j/|\lambda_\mu|^2\rfloor$ so that $Z_\mu^x(t_j) = -\lambda_\mu S_{n_j}(x)$ and $|\lambda_\mu|^2 (n_j-n_{j-1})\to (t_j-t_{j-1})$, and fix $u_1,\dots,u_\ell\in \R$.
By \eqref{eq:diffusive-block-limit}, we have 
\[ \PP_{-iu_j\lambda_\mu,\mu}^{\,n_j-n_{j-1}} \to e^{-(t_j - t_{j-1}) (iu_j+\frac{\sigma_0^2}{2} u_j^2)} \Pi_{0,\mu_0} \quad\text{in operator norm}.  \]
By successively conditioning on the sigma-algebras $\mathcal F_{n_j}$, this implies that the characteristic function of the increments of $Z_\mu^x$ satisfies
\begin{align}
&\E_\mu\exp\!\left(i\sum_{j=1}^\ell
u_j[Z_\mu^x(t_j)-Z_\mu^x(t_{j-1})]\right)\notag\\
&\quad=\E_\mu\exp\!\left(-i\sum_{j=1}^\ell
u_j\lambda_\mu(S_{n_j}(x)-S_{n_{j-1}}(x))\right)\notag\\
&\quad=\bigl(\PP_{-iu_1\lambda_\mu,\mu}^{\,n_1}\cdots
\PP_{-iu_\ell\lambda_\mu,\mu}^{\,n_\ell-n_{\ell-1}}
\mathbbm{1}_{\mathbf{P}(\mathbb{R}^d)}\bigr)([x])\notag\\
&\quad\longrightarrow\prod_{j=1}^\ell
 e^{(t_j-t_{j-1})(-iu_j-\sigma_0^2u_j^2/2)}.
\label{eq:diffusive-increment-transform}
\end{align}
Note that the limit is uniform in $x\in \mathbb{S}^{d-1}$. The limit coincides with the joint characteristic function of the increments of $(Y_{t_j})$ since $Y_{t_j} - Y_{t_{j-1}} \sim \mathcal{N}(-(t_j - t_{j-1}), \sigma_0^2 (t_j - t_{j-1}))$ and the increments are independent (hence the product). This establishes convergence of the finite dimensional distributions.

The second condition requires that $(\sigma_0 B_K-K) - (\sigma_0 B_{K-\delta} - (K-\delta))$ converge in distribution to $0$ as $\delta\to 0$. The increment actually converges almost surely to $0$, which implies convergence in distribution.

We now prove the tightness criterion \eqref{tightness_Cond3}, which verifies the third condition in \cite[Theorem 8.3]{Sagitov2020} for $\alpha=\beta=2$ and $H(x) = 2C_K x$. We start with the following second moment bound: for all $m\geq 1$, $x\in \mathbb{S}^{d-1}$ and nearby noncritical measures $\mu$, we have
\begin{equation}\label{Sec_Moment_Bound}
    \E_\mu [|S_m(x)-m\lambda_\mu|^2] \leq C m,
\end{equation}
for some uniform constant $C$.
Indeed, from the Taylor expansion in \eqref{diffusive-spectral-input}, we know that $\Lambda_\mu(q)-q\lambda_\mu = O(q^2)$ uniformly in $\mu$ for real $q$ with $|q|\le q_0$, where $q_0>0$ is chosen strictly inside the spectral neighbourhood, and $\|\mathcal Q_{q,\mu}^{\,m} \|_{\mathcal C^\alpha}\le C\rho^m$. Choose $\rho_1\in (\rho,1)$. Substituting these estimates into the spectral decomposition, we obtain
\begin{equation}\label{eq:residence-centered-transform}
\begin{aligned}
\E_\mu e^{q(S_m(x)-m\lambda_\mu)}
&=e^{-m\lambda_\mu q}(\PP_{q,\mu}^{\,m}\mathbbm{1}_{\mathbf{P}(\mathbb{R}^d)})([x])\\
&=e^{m(\Lambda_\mu(q)-\lambda_\mu q)}
(\Pi_{q,\mu}\mathbbm{1}_{\mathbf{P}(\mathbb{R}^d)})([x])\\
&\quad+e^{-m\lambda_\mu q}(\mathcal Q_{q,\mu}^{\,m}\mathbbm{1}_{\mathbf{P}(\mathbb{R}^d)})([x])\\
&\le Ce^{C_1mq^2}+Ce^{mq_0|\lambda_\mu|}\rho^m\\
&\le Ce^{C_1mq^2}+C\rho_1^m
\le C_0e^{C_1mq^2},\qquad q\in[-q_0,q_0],
\end{aligned}
\end{equation}
since $\lambda_\mu\to 0$, we may shrink the law neighbourhood so that $\rho e^{q_0|\lambda_\mu|}\le\rho_1$ (the constants do not depend on $x$, $m$ and $\mu$).
Now, for all $y\in \mathbb{R}$, we have $y^2\le e^y + e^{-y}$. Setting $y=\frac{q_0}{\sqrt{m}} (S_m(x)-m\lambda_\mu)$, taking expectations and using the bound above, we obtain
\begin{equation}
\begin{aligned}\label{eq:residence-second-moment}
\frac{q_0^2}{m} \E_\mu|S_m(x)-m\lambda_\mu|^2
&\leq \E_\mu e^{\frac{q_0}{\sqrt{m}}(S_m(x)-m\lambda_\mu)}\\
&\quad+\E_\mu e^{-\frac{q_0}{\sqrt{m}}(S_m(x)-m\lambda_\mu)} \\
& \leq 2 C_0 e^{C_1 q_0^2},
\end{aligned}
\end{equation}
uniformly in $x$ and $\mu$.
This proves \eqref{Sec_Moment_Bound}. We can now prove the tightness criterion. Fix $K>0$ and write $h=|\lambda_\mu|^2$ and $Z = Z_\mu^x$ to simplify notation. Let $(\mathcal{G}_t)$ be the filtration generated by the matrices applied by the rescaled time $t$, i.e. $\mathcal{G}_t = \mathcal{F}_{\lfloor t/h\rfloor}$ ($\mathcal{F}_n = \sigma(A_1,\dots,A_n)$ as in the previous sections). For $0\le s \le t \le K$, let $m$ be the number of matrix steps happening in $(s,t]$. It is easy to see that $mh\le t-s+h$ and $mh\le K$. Set $v_s = A_\omega^{\lfloor s/h \rfloor}x/\|A_\omega^{\lfloor s/h \rfloor}x\|$. Using these simple time bounds, \eqref{Sec_Moment_Bound} and independence of the future matrices, we have 
\begin{equation}\label{eq:residence-conditional-second-moment}
\begin{aligned}
\E_\mu\!\left[|Z(t)-Z(s)|^2\mid\mathcal G_s\right]
&=h\,\left.\E_\mu[S_m(v)^2]\right|_{v=v_s}\\
&=h\,\left.\E_\mu[(S_m(v)-m\lambda_\mu+m\lambda_\mu)^2]\right|_{v=v_s}\\
&\le 2Chm+2(mh)^2\\
&\le C_Kmh\le C_K(t-s+h).
\end{aligned}
\end{equation}
The expectations on the right are taken over a fresh product started at a fixed unit vector $v$, before evaluating at $v=v_s$. Here $C_K=2(C+K)$ is independent of $\mu$ and the initial direction $x$.
Finally, for $0\le r \le s \le t \le K$, conditioning at time $s$ then time $r$ and applying \eqref{eq:residence-conditional-second-moment} twice yields
\begin{equation}\label{eq:residence-mixed-moment}
\begin{aligned}
&\E_\mu\!\left[|Z(s)-Z(r)|^2|Z(t)-Z(s)|^2\right]\\
&\quad=\E_\mu\!\left[|Z(s)-Z(r)|^2
\E_\mu\!\left[|Z(t)-Z(s)|^2\mid\mathcal G_s\right]\right]\\
&\quad\le C_K(t-s+h)\E_\mu|Z(s)-Z(r)|^2\\
& \quad= C_K(t-s+h) \E_\mu \left[   \E_\mu [  |Z(s)-Z(r)|^2 \mid \mathcal{G}_r  ]    \right]  \\
&\quad\le C_K^2(s-r+h)(t-s+h).
\end{aligned}
\end{equation}
If $t-r<h$, at most one matrix step occurs in $(r,t]$, so one of the
two increments is zero. If $t-r\ge h$, each of the factors
$s-r+h$ and $t-s+h$ is at most $2(t-r)$. Hence
\eqref{eq:residence-mixed-moment} gives, in all cases,
\begin{equation}\label{eq:residence-tightness-bound}
\E_\mu\!\left[|Z(s)-Z(r)|^2|Z(t)-Z(s)|^2\right]
\le 4C_K^2(t-r)^2.
\end{equation}
This establishes the third condition, and since all estimates are uniform in $x\in \mathbb{S}^{d-1}$, this proves the path limit along any sequence of noncritical laws $\mu_j$ and starting directions $x_j\in \mathbb{S}^{d-1}$.
\end{proof}
The proof of Theorem \ref{opposite-exponents} is based on applying Theorem \ref{lem:diffusive-cocycle-limit} to transfer path convergence to convergence of the occupation times of $Z_\mu^x$ and $\sigma_0 B_t-t$. This is a function of both processes, so we need a technical lemma that will allow us to apply the continuous mapping theorem.
\begin{lemma}\label{continuous_mapping_theorem_helper}
Let $K>0$, and let $\mathcal{A} \subset D([0,K])$ be the subset of continuous paths defined by 
\[ \mathcal{A} = \{z\in C([0,K]): \Leb (\{t\in[0,K]:z(t)=0\}) =0    \}. \]
Then, the Borel measurable map $g_K: D([0,K])\to \R$ defined by 
\[g_K(z) = \int_0^K  \mathbbm{1}_{ \{z(t)>0\} } dt \]
is continuous at every $z\in\mathcal A$ in the $J_1$ topology.
\end{lemma}
\begin{proof}
    Let $z\in \mathcal{A}$ and let $z_n\to z$ in the $J_1$ topology. By \cite[Lemma 6.27]{Sagitov2020}, since $z$ is continuous, $J_1$ convergence implies uniform convergence, so that \[ \sup_{0\le t \le K} |z_n(t)-z(t)|\to 0. \] This implies that for all $t\in [0,K]$ such that $z(t)\ne 0$, we have 
    \[ \mathbbm{1}_{ \{ z_n(t)>0 \}  } \to \mathbbm{1}_{ \{ z(t)>0 \}  }. \]
    But since $z(t)\ne 0$ Lebesgue-almost everywhere, the above holds for Lebesgue-almost every $t\in [0,K]$. By the dominated convergence theorem, we obtain $g_K(z_n) \to g_K(z)$ which proves continuity.
\end{proof}

\begin{proof}[Proof of Theorem \ref{opposite-exponents}]
Write $Z(t) = \sigma_0 B_t -t$. For $K>0$, define the truncated occupation times by time $K$
\[
\TT_{\mu,K}^x:=\int_0^K\mathbbm{1}_{\{Z_\mu^x(t)>0\}}\,dt = g_K(Z_\mu^x),
\qquad
\TT_K:=\int_0^K\mathbbm{1}_{\{\sigma_0B_t-t>0\}}\,dt=g_K(Z),
\]
and the full occupation times
\[
\TT_{\mu}^x:=\int_0^\infty\mathbbm{1}_{\{Z_\mu^x(t)>0\}}\,dt = |\lambda_\mu|^2 \Topp,
\qquad
\TT:=\int_0^\infty\mathbbm{1}_{\{\sigma_0B_t-t>0\}}\,dt.
\]
By Theorem \ref{lem:diffusive-cocycle-limit}, $Z_\mu^x\xrightarrow{\ d\ }Z$ in $D([0,K])$ with the $J_1$ topology.
The paths of $Z(t)$ are almost surely continuous, and \[ \E[ \Leb(\{t\in[0,K]:Z(t)=0\}) ] = \int_0^K \Prob[Z(t)=0] dt = 0, \]
by Fubini's theorem and the fact that $Z(t)$ has a nondegenerate normal distribution for every $t>0$. By Lemma \ref{continuous_mapping_theorem_helper}, we know that $g_K$ is continuous at every $z\in \mathcal{A}$, and the above shows that $\Prob[Z\in \mathcal{A}]=1$. Hence, the probability that $Z$ takes values in the set of discontinuities of $g_K$ is zero. By the continuous mapping theorem applied to $g_K$, we therefore have $g_K(Z_\mu^x)\xrightarrow{\ d\ }g_K(Z)$. We thus obtain convergence in distribution of the truncated occupation times and their expected values (since both truncated occupation times take values in $[0,K]$):
\begin{equation}\label{eq:finite-residence-limit}
\TT_{\mu,K}^x\xrightarrow{\ d\ }\TT_K,
\qquad
\E_\mu[\TT_{\mu,K}^x]\longrightarrow\E[\TT_K].
\end{equation}

We now require some (uniform in $\mu$) tail estimates to pass $K\to\infty$ and to establish that $\E[\TT]$ and $\E_\mu[\TT_\mu^x]$ are finite.

Tail bounds for $\TT_\mu^x:$ for $q_\mu=-\lambda_\mu/\sigma_0^2$ as defined before, $-\operatorname{sgn}(\lambda_\mu) S_n(x)>0$ implies $e^{q_\mu S_n}>1$, so that
\begin{equation}\label{eq:opposite-time-tail}
\Prob_\mu(-\operatorname{sgn}(\lambda_\mu) S_n(x)>0)
\le\E_\mu[e^{q_\mu S_n(x)}]
\le Ce^{-c|\lambda_\mu|^2n}
\end{equation}
by \eqref{eq:weight-negative-drift} and Markov's inequality, and summing over $n$ yields $\E_\mu[\Topp]<\infty$ for all nearby noncritical measures $\mu$. For $n=\lfloor t/|\lambda_\mu|^2\rfloor$, since $|\lambda_\mu|^2 n \ge t-|\lambda_\mu|^2 $, restricting to all measures with $|\lambda_\mu|<1$ gives for all $t>0$:
\[
\Prob_\mu(Z_\mu^x(t)>0)
\le Ce^{-c|\lambda_\mu|^2n}\le Ce^{c|\lambda_\mu|^2}e^{-ct}\le C'e^{-ct}.
\]
By taking expectations and by Tonelli's theorem, this implies
\begin{align}
\E_\mu[\TT_\mu^x-\TT_{\mu,K}^x]&= \E_\mu \left[  \int_K^\infty \mathbbm{1}_{ \{ Z_\mu^x(t)>0 \} } dt  \right]
\\&=\int_K^\infty\Prob_\mu(Z_\mu^x(t)>0)\,dt\notag\\
&\le C'\int_K^\infty e^{-ct}\,dt
=C_{\mathrm{tail}}e^{-cK},
\label{eq:residence-tail-expectation}
\end{align}
where $C_{\mathrm{tail}}$ is independent of $x,K$ and $\mu$.

Tail bounds for $\TT$: we start by explicitly computing $\E[\TT]$.
Write $\Phi$ and $\phi$ for the standard normal distribution function
and density. Tonelli's theorem and the substitution $t=\sigma_0^2u^2$
give
\begin{align}
\E[\TT]
&=\int_0^\infty\Prob(\sigma_0B_t-t>0)\,dt\notag\\
&=\int_0^\infty\Phi\!\left(-\frac{\sqrt t}{\sigma_0}\right)\,dt
=2\sigma_0^2\int_0^\infty u\Phi(-u)\,du\notag\\
&=\sigma_0^2\int_0^\infty z^2\phi(z)\,dz
=\frac{\sigma_0^2}{2}.
\label{eq:brownian-residence-mean}
\end{align}
In particular, $\TT$ is finite almost surely. By the exponential Markov's inequality,
\[
\Prob(\sigma_0B_t>t)
\le e^{-t/\sigma_0^2}\E[e^{B_t/\sigma_0}]
=e^{-t/(2\sigma_0^2)},\qquad t>0.
\]
Consequently,
\begin{equation}\label{eq:brownian-residence-tail}
\E[\TT-\TT_K]
=\int_K^\infty\Prob(\sigma_0B_t-t>0)\,dt
\le 2\sigma_0^2e^{-K/(2\sigma_0^2)}.
\end{equation}
We can finally prove the convergence in distribution and in means of Theorem \ref{opposite-exponents} using an '$\varepsilon/3$' argument. It suffices to test against bounded Lipschitz functions on $[0,\infty)$, so fix such a function $f$ with Lipschitz constant $L$ (i.e. $|f(x)-f(y)|\le L |x-y|$). The triangle inequality gives
\begin{equation}\label{eq:residence-truncation-comparison}
\begin{aligned}
\left|\E_\mu f(\TT_\mu^x)-\E f(\TT)\right|
&\le\left|\E_\mu f(\TT_\mu^x)-\E_\mu f(\TT_{\mu,K}^x)\right|\\
&\quad+\left|\E_\mu f(\TT_{\mu,K}^x)-\E f(\TT_K)\right|\\
&\quad+\left|\E f(\TT_K)-\E f(\TT)\right|\\
&\le L\E_\mu[\TT_\mu^x-\TT_{\mu,K}^x]\\
&\quad+\left|\E_\mu f(\TT_{\mu,K}^x)-\E f(\TT_K)\right|\\
&\quad+L\E[\TT-\TT_K].
\end{aligned}
\end{equation}
Fixing $K>0$ and taking $\|\mu-\mu_0\|_{TV} \to 0$, the second term tends to zero by \eqref{eq:finite-residence-limit}.
The first and last terms are controlled by
\eqref{eq:residence-tail-expectation} and
\eqref{eq:brownian-residence-tail}. Therefore, first letting
$\mu\to\mu_0$ with $K$ fixed, we obtain
\begin{equation}\label{eq:residence-truncation-limsup}
\begin{aligned}
\limsup_{\substack{\mu\to\mu_0\\\lambda_\mu\ne0}}
\left|\E_\mu f(\TT_\mu^x)-\E f(\TT)\right|
&\le LC_{\mathrm{tail}}e^{-cK}
+2L\sigma_0^2e^{-K/(2\sigma_0^2)}.
\end{aligned}
\end{equation}
Now let $K\to\infty$. The right-hand side vanishes, proving $\E_\mu f(\TT_\mu^x)\to\E f(\TT)$ for every bounded Lipschitz $f$ uniformly in $x\in \mathbb{S}^{d-1}$. This proves \eqref{eq:opposite-residence-law}.
For the convergence of means, the same argument applies by writing 
\begin{equation}\label{eq:residence-mean-comparison}
\begin{aligned}
\left|\E_\mu[\TT_\mu^x]-\E[\TT]\right|
\le\E_\mu[\TT_\mu^x-\TT_{\mu,K}^x]
+\left|\E_\mu[\TT_{\mu,K}^x]-\E[\TT_K]\right|
+\E[\TT-\TT_K].
\end{aligned}
\end{equation}
We obtain \[
|\lambda_\mu|^2\E_\mu[\Topp]
=\E_\mu[\TT_\mu^x]\longrightarrow\E[\TT]
=\frac{\sigma_0^2}{2},
\]
uniformly in $x\in \mathbb{S}^{d-1}$, where the last equality is \eqref{eq:brownian-residence-mean}. This proves the critical exponent result \eqref{eq:opposite-residence-limit}.

\end{proof}

\bibliographystyle{plain}
\bibliography{ERMP}

\end{document}